\documentclass[11pt]{amsart}
\usepackage{amsmath,amsthm,amssymb,amsfonts}
\usepackage[hidelinks]{hyperref}
\usepackage{cleveref}
\usepackage[margin=1in]{geometry}
\usepackage{mathrsfs}
\usepackage{tikz-cd}
\usepackage{pdflscape}
\usepackage{makecell}
\usepackage{longtable}
\usepackage{xcolor}
\usepackage{multirow}
\usepackage{hyperref}

\newtheorem{theorem}{Theorem}[section]
\newtheorem{lemma}[theorem]{Lemma}
\newtheorem{proposition}[theorem]{Proposition}

\theoremstyle{definition}

\theoremstyle{definition}
\newtheorem{remark}[theorem]{Remark}

\newcommand{\Z}{\mathbb{Z}}
\newcommand{\Q}{\mathbb{Q}}

\newcommand{\PP}{\mathbb{P}}
\newcommand{\OO}{\mathcal{O}}

\title[Del Pezzo orbifolds]{Classification of Low degree del Pezzo orbifolds}
\author{Saptarshi Dandapat}
\address{Graduate School of Mathematics, Nagoya University, Furocho Chikusa-ku, Nagoya, 464-8601, Japan}
\email{saptarshi.dandapat.d2@math.nagoya-u.ac.jp}
\date{\today}

\begin{document}

\maketitle

\begin{abstract}
    In this paper we classify low degree del Pezzo orbifolds with irreducible boundaries. In order to achieve desired boundaries,
    we classify low degree curves on low degree del Pezzo surfaces. The notion of Campana orbifolds was introduced by Campana in 2004. A del Pezzo orbifold is a Campana orbifold whose underlying surface is a del Pezzo surface.
    The classification is elementary applications of adjunction formula, Riemann-Roch theorem, Hodge Index theorem and Kawamata-Viehweg vanishing theorem. 
\end{abstract}

\tableofcontents

\section{Introduction} The existence of rational and integral points on algebraic varieties over the function field of a smooth projective curve over an algebraically closed field has been a central object of interest in arithmetic algebraic geometry.
In recent years, the theory of Campana points has become a focus of attention not only because it provides a connection between rational and integral points, but also due to its geometric and arithmetic perspectives.
By the valuative criterion, this leads to the study of Campana sections of fibrations over curves. 

The theory of Campana curves/sections and the validity of weak approximation in the setting of Campana sections have been well studied by Chen, Lehmann and Tanimoto in \cite{chen2024campana}. They proved that for any fibration whose general fibers satisfy a version of Campana rational connectedness (which is a stronger property, equivalent to the notion of orbifold rational connectedness introduced by Campana, e.g., \cite{Campana2010}, \cite{Campana2011}, \cite{Campana2011Survey}) and weak approximation holds for Campana sections at places of good reduction, and they have verified this property for toric Campana orbifolds with toric boundaries. 
A key idea to prove Campana rational connectedness for a Fano orbifold $(X, \Delta)$ is to show an existence of a free Campana curve $f: C \to X$ such that $f_*[C]$ is in the interior of the nef cone $\operatorname{Nef}_1(X)$ of curves, due to the rational connectedness of Fano varieties.

In this paper, we establish the classification of low degree ($\leq 5$) del Pezzo orbifolds when the boundary is irreducible. Which lays a grundwork to show an existence of such free Campana curve on a low degree del Pezzo surface. We also classify low degree weak del Pezzo orbifolds in this setup.

\subsection{\textbf{Campana orbifolds.}} Throughout this paper, we work over an algebraically closed field $\mathbf{k}$ of arbitrary characteristic $p = \operatorname{char} \mathbf{k}$. For simplicity, we assume that $\mathbf{k}$ has characteristic $0$ in this introduction.

Let $\underline{X}$ be a smooth projective variety with a strict normal crossings (SNC) divisor $\Delta = \cup_i \Delta_i$. Let $X = (\underline{X}, \mathcal{M}_X)$ be the log scheme associated to the pair $(\underline{X}, \Delta)$, where $\mathcal{M}_X$ is the sheaf of monoids on $\underline{X}$ defined by 
$$\mathcal{M}_X(U) := \{f \in \mathcal{O}_{\underline{X}}(U) \mid f|_{U \setminus \Delta} \in \mathcal{O}_{\underline{X}}^{\times}(U \setminus \Delta)\}.$$ For each irreducible component $\Delta_i$ of $\Delta$, we assign a weight $\epsilon_i = 1 - \frac{1}{m_i}$, where $m_i \geq 1$ is an integer. We then define the effective $\Q$-divisor $$\Delta_{\epsilon} = \sum\limits_i \epsilon_i \Delta_i.$$
The pair $(X, \Delta_{\epsilon})$ is called a klt Campana orbifold (in the sense of \cite{campana2004orbifolds}).
For a Campana orbifold we always assume that the pair $(X, \Delta)$ is log smooth, and we use the adjective ``klt'' to emphasize that the coefficients of $\Delta_{\epsilon}$ are less than $1$.
In this paper we discuss the case when $\Delta_{\epsilon}$ has only one component, and we use the notation $\Delta_{\epsilon} = \epsilon D$, where $D$ is the only irreducible component, which is a $\Q$-divisor on $\underline{X}$. 

\begin{remark}
    The equivalence between the above notions of ``orbifold'' and the notion due to Campana is discussed in \cite[Remark~7.3]{chen2024campana}.
\end{remark}

\textbf{Del Pezzo orbifolds.} Let $\underline{X}$ be a smooth del Pezzo surface. A del Pezzo orbifold is a Campana orbifold $(X, \Delta_{\epsilon})$ such that $-(K_{\underline{X}} + \Delta_{\epsilon})$ is ample.

In this paper we classify all del Pezzo orbifolds with irreducible boundary, that is the case when $\Delta_{\epsilon} = \epsilon D$, for an irreducible $\Q$-divisor $D$ on $\underline{X}$. In particular, we are interested in the case when $-(K_{\underline{X}} + \epsilon D)$ is ample.

\subsection{\textbf{Main results.}} Assume that $\mathbf{k}$ has characteristic $0$. Let $(X, \Delta_{\epsilon})$ be a klt del Pezzo orbifold over $\mathbf{k}$. For the rest of the paper we assume $\Delta_{\epsilon} = \epsilon D$ where $D$ is an irreducible $\Q$-divisor.

\begin{theorem}
    Let $(X, \epsilon D)$ is a klt (weak) del Pezzo orbifold such that $\underline{X}$ has degree $d \leq 5$ and $D$ is irreducible $\Q$-divisor of anticanonical degree $-K_{\underline{X}} \cdot D \leq 2d$. Then the pair $(X, \epsilon D)$ must be one of the following options.
    \begin{enumerate}
        \item{Section 3.} $d = 1$. Then $-(K_{\underline{X}} + \epsilon D)$ is ample for all $0 < \epsilon < 1$ if $D$ is a member of the anticanonical class $|-K_{\underline{X}}|$. It is nef but not big if $D$ is a member of $|-2K_{\underline{X}}|$ and $\epsilon = \frac{1}{2}$.
        \item{Section 4.} $d = 2$. Then $-(K_{\underline{X}} + \epsilon D)$ is ample for all $0 < \epsilon < 1$ if $D$ is a member of the anticanonical class $|-K_{\underline{X}}|$. It is big and nef if $D$ is a $(-1)$-curve on $\underline{X}$ and $\epsilon = \frac{1}{2}$. It is nef but not big if $D$ is either a smooth fiber of a conic fibration or is a member of $|-2K_{\underline{X}}|$ and $\epsilon = \frac{1}{2}$.
        \item{Section 5.} $d = 3$. Then $-(K_{\underline{X}} + \epsilon D)$ is ample for all $0 < \epsilon < 1$ if $D$ is either a $(-1)$-curve, or a member of the anticanonical class $|-K_{\underline{X}}|$. It is big and nef if, either $D$ pullback of the hyperplane class in $\PP^2$, or $D$ is smooth fiber of a conic fibration and $\epsilon = \frac{1}{2}$.
        It is nef but not big if $\epsilon = \frac{1}{2}$ and $D$ is either pullback of the anticanonical divisor on a quartic del Pezzo surface or is a member of $|-2K_{\underline{X}}|$.
        \item{Section 6.} $d = 4$. Then $-(K_{\underline{X}} + \epsilon D)$ is ample for all $0 < \epsilon < 1$ if $D$ is either a $(-1)$-curve, or a smooth fiber of a conic fibration, or a member of the anticanonical class $|-K_{\underline{X}}|$. It is big and nef if $\epsilon = \frac{1}{2}$ and $D$ is either pullback of a hyperplane class in $\PP^2$, or pullback of a member of $(1,1)$ class in $\PP^1 \times \PP^1$. 
        It is nef but not big if $\epsilon = \frac{1}{2}$ and $D$ is either linearly equivalent to $-K_{\underline{X}} + F$ for some smooth conic $F$, or is a member of $|-2K_{\underline{X}}|$.
        \item{Section 7.} $d = 5$. Then $-(K_{\underline{X}} + \epsilon D)$ is ample for all $0 < \epsilon < 1$ if $D$ is either a $(-1)$-curve, or a smooth fiber of a conic fibration, or pullback of a line in $\PP^2$, or a member of the anticanonical class. 
        It is big and nef if $\epsilon = \frac{1}{2}$ and $D$ is either pullback of a member in the $(1,1)$-class of $\PP^1 \times \PP^1$, or is pullback of a divisor of the form $2H' - E$ from blow up of $\PP^2$ at one point, where $H'$ is the pullback of the hyperplane class in $\PP^2$ via a blow up $\beta':\underline{X}' \to \PP^2$ at one point and $E$ is the exceptional divisor of $\beta'$, or pullback of a member in the anticanonical class of a smooth sextic del Pezzo surface, or is linearly equivalent to $-K_{\underline{X}} + F$ for a smooth conic $F$.
        It is nef but not big if $\epsilon = \frac{1}{2}$ and $D$ is pullback of a divisor in the $(1,2)$ or $(2,1)$-class in $\PP^1 \times \PP^1$ via a birational morphism to $\PP^1 \times \PP^1$. In particular, they are of the form:
        $2H$ or $4H - 2E_i - 2E_j - 2E_k$ where $H$ is the pullback of a hyperplane class in $\PP^2$ via a blow up $\beta : \underline{X} \to \PP^2$ at four general points, as described in section \ref{dp5}, and $E_i, E_j, E_k, \{0 \leq i,j,k \leq 4\}$ are exceptional divisors of $\beta$, or $D + K_{\underline{X}}$ is pullback of a line in $\PP^2$, or $D$ is a member of $|-2K_{\underline{X}}|$.
    \end{enumerate}  
\end{theorem}

At the end of the paper, we record the classification of irreducible curves of anticanonical degree less than or equal to $2d$ by their self intersection numbers, on a del Pezzo surface of degree $d \leq 5$. 

\vspace{0.5cm}

\textbf{Acknowledgements:} The author is grateful to his advisor, Sho Tanimoto, for introducing the problem. He also thanks Sho Tanimoto, Runxuan Gao, Marta Pieropan, Brian Lehmann, and Qile Chen for valuable comments and discussions on the first draft of the paper. The author further thanks Sho Tanimoto, Qile Chen, and Brian Lehmann for sharing their unpublished work on the problem in the case of cubic surfaces.

The author was supported by a Japanese Government (MEXT) Scholarship for Research Students (Doctoral Course), recommended by the Embassy of Japan in India (New Delhi).

\section{Preliminaries}

\subsection{} We recall the following properties of a del Pezzo surface over $\mathbf{k}$. A del Pezzo surface is a smooth projective rational surface $\underline{X}$ such that the anticanonical class $-K_{\underline{X}}$ is ample. If $\underline{X}$ is over an algebraically closed field then $\underline{X}$ is either one of the following: $\PP^2, \PP^1 \times \PP^1$ or $\PP^2$ blown up at up to $8$ general points. Points are said to be in general position if no three are collinear, no six lie on a conic, and no eight lie on a cubic with one of them a singular point of the cubic. 
A del Pezzo surface of degree $d$ can be obtained by blowing up $\PP^2$ at $9-d$ general points. 

\begin{enumerate}
    \item{Del Pezzo surface of degree 1.} It is isomorphic to the blow up of $\PP^2$ at eight general points. If the characteristic of $k$ is not $2$ or $3$, then it can be given by a hypersurface of degree $6$ in the weighted projective space $\PP(1,1,2,3)$ given by $w^2 = t^3 + f_4(x,y)t +f_6(x,y)$, with $\operatorname{deg}(f_i) = i$. In this case the anticanonical divisor has degree $1$.
    There are $240$ $(-1)$-curves in this surface.
    \item{Del Pezzo surface of degree 2.} It is a double cover of $\PP^2$ ramified over a smooth quartic curve $Q$, (this notion also works in characteristic $2$ by \cite[Proposition 3.1]{dolgachev2025automorphisms}) which is isomorphic to the blow up of $\PP^2$ at seven general points. In this case the anticanonical divisor has degree $2$.
    There are $56$ $(-1)$-curves in this surface.
    \item{Del Pezzo surface of degree 3.} It is a smooth cubic surface in $\PP^3$, which is isomorphic to the blow up of $\PP^2$ at six general points. In this case the anticanonical divisor has degree $3$. 
    There are $27$ $(-1)$-curves in this surface.
    \item{Del Pezzo surface of degree 4.} It is isomorphic to the blow up of $\PP^2$ at five general points. In this case the anticanonical divisor has degree $4$.
    There are $16$ $(-1)$-curves in this surface.
    \item{Del Pezzo surface of degree 5.} It is isomorphic to the blow up of $\PP^2$ at four general points. In this case the anticanonical divisor has degree $5$.
    There are $10$ $(-1)$-curves in this surface.
\end{enumerate}

\subsection{} For del Pezzo surfaces of degree $\leq 7$, we record the following lemma from \cite[Proposition 3.9]{derenthal2008nef} for our computations.

\begin{lemma}
    \label{thelemma}
    Let $\underline{X}$ be a smooth del Pezzo surface of degree $\leq 7$. Then $\operatorname{Eff}_1(\underline{X})$ is generated by $(-1)$-curves and $\operatorname{Nef}^1(\underline{X})$ is generated by smooth conics and the pullbacks of the hyperplane class via a birational morphism to $\PP^2$.
\end{lemma}

\subsection{} To classify the irreducible curves on each del Pezzo surface, we find a suitable birational morphism to another del Pezzo surface, depending on its anticanonical degree and self-intersection number, such that it can be realised as pullback of some known curve. 
Note that if a class is represented by an irreducible curve, then it is nef unless it is represented by a $(-1)$-curve, and a nef divisor $D$ is big if and only if its top self-intersection is strictly positive.

Since a del Pezzo surface is a rational surface, we have $\chi(\OO_X) = 1$ for any del Pezzo surface $X$. Using Kawamata-Viehweg vanishing and Riemann-Roch we determine whether a particular curve $D$ on a surface $X$ is effective, that is $h^0(X, \OO_X(D)) > 0$. 
By \cite[Theorem~2.3.19]{lazarsfeld2017positivity} if $D$ is a pseudoeffective integral divisor we on $X$ then by Zariski decomposition $D$ can be uniquely written as a sum $$D = F + \sum a_i N_i, \quad a_i > 0 \text{ for all } i$$ of $\Q$-divisors with the following properties:
\begin{enumerate}
    \item $F$ is nef.
    \item $N = \sum\limits_i a_i N_i$ is effective, and if $N \neq 0$ then the intersection matrix $$\lVert (N_i \cdot N_j) \rVert$$ determined by the components of $N$ is negative definite.
    \item $F$ is orthogonal to each of the components of $N$, i.e. $(F \cdot N_i) = 0$ for every $i$.
\end{enumerate}

Therefore we can conclude that $D \cdot N_i < 0$ for all $i$ and $N_i$'s are disjoint $(-1)$-curves.

\subsection{}\label{prog} For a del Pezzo surface $\underline{X}$ of degree $d$, we find the explicit form of an irreducible curve $C$ with anticanonical degree $m$ and $C^2 = n$ using a Python 3 program using the technique of \cite[V, Theorem 4.9]{hartshorne2013algebraic}. 
If $-K_{\underline{X}} \cdot C = m$ and $C^2 = n$, we can write $C \sim aH - \sum\limits_{i = 1}^{9-d} b_i E_i$ where $H$ is pullback of a hyperplane class in $\PP^2$ via a birational morphism $\beta: \underline{X} \to \PP^2$ and $E_i$'s are the exceptional divisors of $\beta$.
Then we have 
\begin{equation*}
    \begin{aligned}
        3a - \sum b_i &= m \\
        a^2 - \sum b_i^2 &= n
    \end{aligned}
\end{equation*}
Using Cauchy-Schwarz's inequality, we get $$\left( \sum b_i \right)^2 \leq (9-d) \left( \sum b_i^2 \right).$$
Substituting $\sum b_i = 3a - m$ and $\sum b_i^2 = a^2 - n$, we get a quadratic equation in $a$, from which we can find the possible integer values of $a$. Hence we find the explicit forms of the curve $C$.
For such a curve, we check the intersection numbers with different $(-1)$-curves by another Python 3 program, to verify our results obtained from geometric configuration of $(-1)$-curves on $\underline{X}$. 

\section{Surfaces of degree 1.}
\label{dp1}

A del Pezzo surface $\underline{X}$, of degree $1$ is isomorphic to the blow up $\beta: \underline{X} \to \PP^2$ at eight general points in $\PP^2$.
Then the canonical divisor of $\underline{X}$ is given by $$K_{\underline{X}} = -3H + E_1 + E_2 + E_3 + E_4 + E_5 + E_6 + E_7 +E_8,$$ where $E_i$'s $(1 \leq i \leq 8)$ are the exceptional curves under $\beta$ and $\beta^*K_{\PP^2} = -3H$, where $H$ is the hyperplane class in $\underline{X}$.
We classify possible del Pezzo orbifolds $(\underline{X}, \Delta_{\epsilon})$ where $\underline{X}$ is as above and $\Delta_{\epsilon}$ is a $\Q$-divisor. In this case, the anticanonical divisor $-K_{\underline{X}}$ has degree $1$, $\Delta_{\epsilon}$ can have at most $2$ components and each component has anticanonical degree at most $2$.

\subsection{When boundary is irreducible.} We classify the case when $\Delta_{\epsilon}$ has only one component first. That is we classify curve classes in $\underline{X}$ which has anticanonical degree at most $2$.

\begin{itemize}
    \item An irreducible curve $C$, of anticanonical degree $1$ is either a $(-1)$-curve satisfying $C^2 = -1,$ or is member of the anticanonical class $|-K_{\underline{X}}|$. In the former case, there are $240$ such curves of the following forms. Using \ref{prog} we get:
    \begin{enumerate}
        \item $E_i$ for $i = 1, \ldots, 8$ are $8$ exceptional curves on $\underline{X}$ under the blow up $\beta: \underline{X} \to \PP^2$ at eight general points in $\PP^2$.
        \item $H - E_i - E_j$, where $1 \leq i, j \leq 8, i \neq j$. There are $28$ such curves.
        \item $2H - E_i - E_j - E_k - E_l - E_m$, where $1 \leq i, j, k, l, m \leq 8$ are distinct indices. There are $56$ such curves.
        \item $3H - 2E_i - E_j - E_k - E_l - E_m - E_n - E_p$. There are $56$ such curves.
        \item $4H - 2E_i - 2E_j - 2E_k - E_l - E_m - E_n - E_p - E_q$. There are $56$ such curves.
        \item $5H - 2E_i - 2E_j - 2E_k - 2E_l - 2E_m - 2E_n - E_p - E_q$. There are $28$ such curves.
        \item $6H - 3E_i - 2E_j - 2E_k - 2E_l - 2E_m - 2E_n - 2E_p - 2E_q$, where $1 \leq i, j, k, l, m, p, q \leq 8$ are distinct indices. There are $8$ such curves.
    \end{enumerate}
    In the later case, $C$ is be a member of the anticanonical class $|-K_{\underline{X}}|$, in this case $C^2 = 1$.
    \item If $C$ is an irreducible curve of anticanonical degree $2$, then using Hodge Index theorem, $C^2 \leq 4.$ Again using adjunction formula, $2p_a(C) - 2 = C^2 + K_{\underline{X}} \cdot C = C^2 - 2.$ Hence, if $C$ is an effective irreducible curve of anticanonical degree $2$, then $C^2 = 0, 2, 4$.
    \begin{proposition}
        Let $\underline{X}$ be a del Pezzo surface of degree $1$ and $C$ is an irreducible curve of anticanonical degree $2$. Then $C^2 = 0, 2$ or $4$. When $C^2 = 0$, $C$ is smooth fiber of a conic fibration to $\PP^1$. When $C^2 = 2$, $C$ is pullback of a member in the anticanonical class of a del Pezzo surface of degree $2$ via a birational morphism. When $C^2 = 4$, $C$ is a member of $|-2K_{\underline{X}}|$.
    \end{proposition}
    \begin{proof}
        \begin{enumerate}
            \item \textit{Case 1: $C^2 = 0$.} In this case $|C|$ is a basepoint-free pencil. Hence it defines a morphism $f: \underline{X} \to \PP^1$. Such $C$ is a smooth fiber of a conic fibration $f$ on $\underline{X}$. They are one of the following forms.
            Using \ref{prog} we get:
                \begin{enumerate}
                    \item $H - E_i$, where $E_i, 1 \leq i \leq 8$.
                    \item $2H - E_i - E_j - E_k - E_l$, where $1 \leq i, j, k, l \leq 8$ are distinct indices.
                    \item $3H - 2E_i - E_j - E_k - E_l - E_m - E_n$, where $1 \leq i, j, k, l, m, n \leq 8$ are distinct indices.
                    \item $4H - 2E_i - 2E_j - 2E_k - E_l - E_m - E_n - E_p;$ $4H - 3E_i - E_j - E_k - E_l - E_m - E_n - E_p - E_q$.
                    \item $5H - 2E_i - 2E_j - 2E_k - 2E_l - 2E_m - 2E_n - E_p;$ $5H - 3E_i - 2E_j - 2E_k - 2E_l - E_m - E_n - E_p - E_q$.
                    \item $6H - 3E_i - 3E_j - 2E_k - 2E_l - 2E_m - 2E_n - E_p - E_q$.
                    \item $7H - 3E_i - 3E_j - 3E_k - 3E_l - 2E_m - 2E_n - 2E_p - E_q;$ $7H - 4E_i - 3E_j - 2E_k - 2E_l - 2E_m - 2E_n - 2E_p - 2E_q$.
                    \item $8H - 3E_i - 3E_j - 3E_k - 3E_l - 3E_m - 3E_n - 3E_p - E_q;$ $8H - 4E_i - 3E_j - 3E_k - 3E_l - 3E_m - 2E_n - 2E_p - 2E_q$.
                    \item $9H - 4E_i - 4E_j - 3E_k - 3E_l - 3E_m - 3E_n - 3E_p - 2E_q$.
                    \item $10H - 4E_i - 4E_j - 4E_k - 4E_l - 3E_m - 3E_n - 3E_p - 3E_q$.
                    \item $11H - 4E_i - 4E_j - 4E_k - 4E_l - 4E_m - 4E_n - 4E_p - 3E_q$, where $1 \leq i, j, k, l, m, n, p, q \leq 8$ are distinct indices.
                \end{enumerate}
            \item \textit{Case 2: $C^2 = 2$.} In this case $C+K_{\underline{X}}$ is linearly equivalent to a $(-1)$-curve. Let $\pi: \underline{X} \to \underline{X}'$ be a contraction of this $(-1)$-curve. Then $\underline{X}'$ is a del Pezzo surface of degree $2$ and $C$ is a member of $|-\pi^*K_{\underline{X}'}|$. By contracting different $(-1)$-curves, we get one of the following forms.
            Using \ref{prog} we get:
                \begin{enumerate}
                    \item $3H - E_i - E_j - E_k - E_l - E_m - E_n - E_p$, where $1 \leq i, j, k, l, m, n, p \leq 8$ are distinct indices.
                    \item $4H - 2E_i - 2E_j - E_k - E_l - E_m - E_n - E_p - E_q$.
                    \item $5H - 2E_i - 2E_j - 2E_k - 2E_l - 2E_m - E_n - E_p - E_q$.
                    \item $6H - 3E_i - 2E_j - 2E_k - 2E_l - 2E_m - 2E_n - 2E_p - E_q$.
                    \item $7H - 3E_i - 3E_j - 3E_k - 2E_l - 2E_m - 2E_n - 2E_p - 2E_q$.
                    \item $8H - 3E_i - 3E_j - 3E_k - 3E_l - 3E_m - 3E_n - 2E_p - 2E_q$.
                    \item $9H - 4E_i - 3E_j - 3E_k - 3E_l - 3E_m - 3E_n - 3E_p - 3E_q$, where $1 \leq i, j, k, l, m, n, p, q \leq 8$ are distinct indices.
                \end{enumerate}
            \item \textit{Case 3: $C^2 = 4$.} In this case $C \in |-2K_{\underline{X}}|$.
        \end{enumerate}
    \end{proof}
\end{itemize}

\subsection{Classification of irreducible boundaries.} We let $\Delta = \epsilon D$, where $D$ is an irreducible curve on $\underline{X}$, a del Pezzo surface of degree $1$, and $\epsilon = 1 - 1/m.$ We discuss whether $-(K_{\underline{X}}+\Delta)$ satisfies some positivity of divisors when $D$ is an irreducible curve of anticanonical degree at most $2$.

Using Lemma \ref{thelemma} for a del Pezzo Surface of degree $1$, we can conclude that a divisor $L$ is ample (nef) if and only if $L \cdot E > 0$ (resp. $L \cdot E \geq 0$) for any $(-1)$-curve $E$. We use \ref{prog} to calculate the intersection numbers.

\begin{itemize}
    \item \textit{Case of $-K_{\underline{X}} \cdot D = 1:$} We assume that $D$ is a $(-1)$-curve. In this case $D^2 = -1$. Since $D \cdot E = 0, 1, 2$ or $3$ for any $(-1)$-curve $E$ distinct from $D$, we conclude that $-(K_{\underline{X}} + \epsilon D)$ is not nef for any $\epsilon > 0.$

    If $D \in |-K_{\underline{X}}|$ then $D \cdot E = 1$ for any $(-1)$-curve $E$. Hence $-(K_{\underline{X}}+ \epsilon D)$ is ample for all $\epsilon > 0.$
    \item \textit{Case $-K_{\underline{X}} \cdot D = 2:$} If $D$ is a smooth fiber of a conic fibration, then $D \cdot E = 0, 1 ,2, 3$ or $4$ for any $(-1)$-curve $E$. So $(-K_{\underline{X}} + \epsilon D)$ is not nef for any $\epsilon > 0.$

    If $D$ is pullback of a member in the anticanonical class of a degree $2$ del Pezzo surface via a birational morphism, then $D \cdot E = 0, 1, 2, 3$ or $4$ for any $(-1)$-curve $E$. Hence we conclude that $-(K_{\underline{X}} + \epsilon D)$ is not nef for any $\epsilon > 0$.

    If $D \in |-2K_{\underline{X}}|$, then $D \cdot E = 2$. Hence $-(K_{\underline{X}} + \epsilon D)$ is nef for $\epsilon = \frac{1}{2}$, but it is not big.
\end{itemize}

Therefore we classify all del Pezzo orbifolds (underlying surface of degree $1$) by the following proposition.

\begin{proposition}
    Let $(X, \epsilon D)$ be a del Pezzo orbifold whose underlying surface $\underline{X}$ is a del Pezzo surface of degree $1$. Let $D$ be an irreducible curve of degree less than or equal to $2$, and $\epsilon = 1-1/m, m \in \Z_{\geq 2}$. Then $-(K_{\underline{X}} + \epsilon D)$ is ample if and only if $D$ is a member of the anticanonical class, for all $0 < \epsilon < 1.$ Also, $(-K_{\underline{X}} + \epsilon D)$ is nef only when $D$ is a member of $|-2K_{\underline{X}}|$ and $\epsilon = \frac{1}{2}.$ In this case $-(K_{\underline{X}} + \epsilon D)$ is nef but not big.
\end{proposition}

\section{Surfaces of degree 2.}
\label{dp2}

A del Pezzo Surface of degree $2$ can be realised as a double cover of $\PP^2$ ramified over a smooth quartic curve $Q$. This notion works in characteristic $2$ as well due to \cite[Proposition 3.1]{dolgachev2025automorphisms}. Let $\underline{X}$ be such a del Pezzo Surface. Then it is the blow up $\beta: \underline{X} \to \PP^2$ at seven general points. 
In this case the canonical divisor $K_{\underline{X}}$ is given by $$K_{\underline{X}} = -3H + E_1 + E_2 + E_3 + E_4 + E_5 + E_6 + E_7,$$ where $E_i$'s $(1 \leq i \leq 8)$ are the exceptional curves under $\beta$ and $\beta^*K_{\PP^2} = -3H$, where $H$ is the hyperplane class in $\underline{X}$. 
We classify possible del Pezzo orbifolds $(\underline{X}, \Delta_{\epsilon})$ where $\underline{X}$ is as above and $\Delta_{\epsilon}$ is a $\Q$-divisor. In this case, the anticanonical divisor $-K_{\underline{X}}$ has degree $2$, so $\Delta_{\epsilon}$ can have at most $4$ components and each component has anticanonical degree at most $4$.

\subsection{Low degree curves on $\underline{X}$.} We first classify classes of irreducible curves of anticanonical degree at most $2$.

\begin{itemize}
    \item an irreducible curve $C$ of anticanonical degree $1$ on $\underline{X}$ is a $(-1)$-curve satisfying $C^2=-1$. A del Pezzo surface of degree $2$ can be constructed by blowing up $\PP^2$ at $7$ points. The $(-1)$-curves on a degree $2$ del Pezzo surface are one of the following forms. Using \ref{prog} we get:
    \begin{enumerate}
        \item $E_i$ where $i = 1, \cdots, 7$ are exceptional curves of the blow up $\beta: \underline{X} \to \PP^2$ at $7$ general points in $\PP^2$.
        \item $H - E_i - E_j$, where $1 \leq i, j \leq 7, i \neq j$. There are $21$ such curves.
        \item $2H - E_i - E_j - E_k - E_l - E_m$, where $1 \leq i, j, k, l, m \leq 7$ are distinct indices. There are $21$ such curves.
        \item $3H - 2E_i - E_j - E_k - E_l - E_m - E_n - E_p$, where $1 \leq i, j, k, l, m, n, p \leq 7$ are distinct indices. There are $7$ such curves.
    \end{enumerate}
    \item Let $C$ be an irreducible curve of anticanonical degree $2$. Then by Hodge Index theorem and adjunction formula, it must satisfy either $C^2 = 0$ or $C^2 = 2$.
    \begin{enumerate}
        \item \textit{Case 1: $C^2 = 0.$} In this case $|C|$ is a basepoint-free pencil. Hence it defines a morphism $f: \underline{X} \to \PP^1$. Such a $C$ is a smooth fiber of a conic fibration $f$ on $\underline{X}$. They are one of the following forms. Using \ref{prog} we get:
        \begin{enumerate}
            \item $H - E_i$, where $1 \leq i \leq 7$.
            \item $2H - E_i - E_j - E_k - E_l$, where $1 \leq i, j, k, l \leq 7$ are distinct indices.
            \item $3H - 2E_i - E_j - E_k - E_l - E_m - E_n$, where $1 \leq i, j, k, l, m, n \leq 7$ are distinct indices.
            \item $4H - 2E_i - 2E_j - 2E_k - E_l - E_m - E_n - E_p$.
            \item $5H - 2E_i - 2E_j - 2E_k - 2E_l - 2E_m - 2E_n - E_p$, where $1 \leq i, j, k, l, m, n, p \leq 7$ are distinct indices.
        \end{enumerate}
        \item \textit{Case 2: $C^2 = 2.$} In this case, $C$ is member of $|-K_{\underline{X}}|$.
    \end{enumerate}
\end{itemize}

Now we assume that $C$ is an irreducible curve of anticanonical degree $3$. By Hodge Index theorem and adjunction formula, it must satisfy either $C^2 = 1$ or $C^2 = 3$. 

\begin{proposition}
    \label{dp2d3}
    Let $\underline{X}$ be a del Pezzo surface of degree $2$ and $C$ is an irreducible curve of anticanonical degree $3$. Then $C^2 = 1$ or $3$. In the former case, $C$ is pullback of a line via a birational morphism to $\PP^2$ and in the later case $C$ is the pullback of a member of the anticanonical class on a smooth cubic surface via a birational morphism.
\end{proposition}

\begin{proof}
    \begin{enumerate}
        \item Assume that $C^2 = 1$. Hence $3C$ is big and nef. By Kawamata-Viehweg vanishing and Riemann-Roch $3C + K_{\underline{X}}$ must be effective since $h^0(\underline{X}, \OO_{\underline{X}}(3C + K_{\underline{X}})) = 1 > 0$. We also have 
        $$C \cdot (3C + K_{\underline{X}}) = 0.$$ Thus $3C+K_{\underline{X}}$ must be contracted by a big and semi-ample divisor $C$ so that $3C+K_{\underline{X}}$ is an effective linear combination of disjoint $(-1)$-curves. Let $\pi: \underline{X} \to \underline{X}'$ be the birational contraction induced by $C$. Then $C$ is the pullback of a big and nef divisor $C'$ on $\underline{X}'$ satisfying $3C' + K_{\underline{X}'} \sim 0$.  Since $K_{\underline{X}'}$ is divisible by $3$, we conclude that $\underline{X}'$ is $\PP^2$. Then $C'$ must be a member of the hyperplane class on $\PP^2$.
        For different birational contractions, they are one of the following forms. Using \ref{prog} we get:
        \begin{enumerate}
            \item $H$.
            \item $2H - E_i - E_j - E_k$, where $1 \leq i, j, k \leq 7$ are distinct indices.
            \item $3H - 2E_i - E_j - E_k - E_l - E_m$, where $1 \leq i, j, k, l, m \leq 7$ are distinct indices.
            \item $4H - 2E_i - 2E_j - 2E_k - E_l - E_m - E_n,$ $4H - 3E_i - E_j - E_k - E_l - E_m - E_n - E_p$.
            \item $5H - 2E_i - 2E_j - 2E_k - 2E_l - 2E_m - 2E_n$; $5H - 3E_i - 2E_j - 2E_k - 2E_l - E_m - E_n - E_p$.
            \item $6H - 3E_i - 3E_j - 2E_k - 2E_l - 2E_m - 2E_n - E_p$.
            \item $7H - 3E_i - 3E_j - 3E_k - 3E_l - 2E_m - 2E_n - 2E_p,$ where $1 \leq i, j, k, l, m, n, p \leq 7$ are distinct indices.
            \item $8H - \sum\limits_{i=1}^7 3E_i$.
        \end{enumerate}
        \item Assume that $C^2 = 3$. In this case $C+K_{\underline{X}}$ is linearly equivalent to a $(-1)$-curve. Let $\pi: \underline{X} \to \underline{X}'$ be a birational contraction of this $(-1)$-curve. Then $\underline{X}'$ is a smooth cubic surface and $C$ is a member of $|-\pi^*K_{\underline{X}'}|.$
        For different birational contractions, they are one of the following forms. Using \ref{prog} we get:
        \begin{enumerate}
            \item $3H - E_i - E_j - E_k - E_l - E_m - E_n$, where $1 \leq i, j, k, l, m, n \leq 7$ are distinct indices.
            \item $4H - 2E_i - 2E_j - E_k - E_l - E_m - E_n - E_p$.
            \item $5H - 2E_i - 2E_j - 2E_k - 2E_l - 2E_m - E_n - E_p$.
            \item $6H - 3E_i - 2E_j - 2E_k - 2E_l - 2E_m - 2E_n - 2E_p$, where $1 \leq i, j, k, l, m, n, p \leq 7$ are distinct indices.
        \end{enumerate}
    \end{enumerate}
\end{proof}

Next let us assume that $C$ be an irreducible curve of anticanonical degree $4$ in $\underline{X}$. Using Hodge Index Theorem and adjunction formula we get that $2 \leq C^2 \leq 8.$ Hence possible values of $C^2$ are $2,4,6$ and $8$.

\begin{proposition}
    \label{dp2d4}
    Let $\underline{X}$ be a degree $2$ del Pezzo surface and $C$ is an irreducible curve of anticanonical degree $4$. Then we have either $C^2 = 2, 4, 6$ or $8$.
    When $C^2 = 2$, there exists a birational map $\pi: \underline{X} \to \PP^1 \times \PP^1$ such that $C$ is pullback of a member of a curve in $(1,1)$ class.
    When $C^2 = 4$, there exists a birational map $\pi: \underline{X} \to \underline{X}'$ to the quartic del Pezzo surface such that $C$ is a member of $|-\pi^*K_{\underline{X}'}|$.
    When $C^2 = 6$, then $C$ is ample and one can find a smooth conic $F$ such that $C \sim -K_{\underline{X}} + F$.
    When $C^2 = 8$, $C \in |-2K_{\underline{X}}|$.
\end{proposition}

\begin{proof}
    \begin{enumerate}
        \item If $C^2 = 2$, then $2C$ is big and nef. By Kawamata-Viehweg vanishing and Riemann-Roch $2C+K_{\underline{X}}$ must be effective since $h^0(\underline{X}, \OO_{\underline{X}}(2C+K_{\underline{X}})) = 1 > 0$. We also have $$C \cdot (2C + K_{\underline{X}}) = 0.$$ 
        Thus $2C + K_{\underline{X}}$ must be contracted by a big and semi-ample divisor $C$ such that $2C + K_{\underline{X}}$ is an effective linear combination of disjoint $(-1)$-curves. Let $\pi: \underline{X} \to \underline{X}'$ be the birational contraction induced by $C$. Then $C$ is pullback of a big and nef divisor $C'$ on $\underline{X}'$ satisfying $2C' + K_{\underline{X}'} \sim 0$. Since $K_{\underline{X}'}$ is divisible by $2$, we conclude that $\underline{X}'$ is a smooth quadric surface. Then $C'$ must be a member of the hyperplane class on $\underline{X}'$.
        For different pullbacks, they can be described explicitly by one of the following forms. Using \ref{prog} we get:
        \begin{enumerate}
            \item $2H - E_i - E_j$, where $1 \leq i, j \leq 7, i \neq j$.
            \item $3H - 2E_i - E_j - E_k - E_l$, where $1 \leq i, j, k, l \leq 7$ are distinct indices.
            \item $4H - 2E_i - 2E_j - 2E_k - E_l - E_m,$ $4H - 3E_i - E_j - E_k - E_l - E_m - E_n$.
            \item $5H - 3E_i - 2E_j - 2E_k - 2E_l - E_m - E_n$, where $1 \leq i, j, k, l, m, n \leq 7$ are distinct indices.
            \item $6H - 3E_i - 3E_j - 2E_k - 2E_l - 2E_m - 2E_n,$ $6H - 3E_i - 3E_j - 3E_k - 2E_l - E_m - E_n - E_p,$ $6H - 4E_i - 2E_j - 2E_k - 2E_l - 2E_m - E_n - E_p$.
            \item $7H - 4E_i - 3E_j - 3E_k - 2E_l - 2E_m - 2E_n - E_p$.
            \item $8H - 4E_i - 3E_j - 3E_k - 3E_l - 3E_m - 3E_n - E_p,$ $8H - 4E_i - 4E_j - 3E_k - 3E_l - 2E_m - 2E_n - 2E_p$.
            \item $9H - 4E_i - 4E_j - 4E_k - 3E_l - 3E_m - 3E_n - 2E_p$.
            \item $10H - 4E_i - 4E_j - 4E_k - 4E_l - 4E_m - 3E_n - 3E_p$, where $1 \leq i, j, k, l, m, n, p \leq 7$ are distinct indices.
        \end{enumerate}
        \item If $C^2 = 4$, then $C+K_{\underline{X}}$ is disjoint union of two $(-1)$-curves. Let $\pi: \underline{X} \to \underline{X}'$ be the birational contraction of these $(-1)$-curves. Then $\underline{X}'$ is a degree $4$ del Pezzo surface and $C$ is a member of $|-\pi^*K_{\underline{X}'}|$. By contracting different $(-1)$-curves, we get the following forms of $C$. Using \ref{prog} we get:
        \begin{enumerate}
            \item $3H - E_i - E_j - E_k - E_l - E_m$, where $1 \leq i, j, k, l, m \leq 7$ are distinct indices.
            \item $4H - 2E_i - 2E_j - E_k - E_l - E_m - E_n$, where $1 \leq i, j, k, l, m, n \leq 7$ are distinct indices.
            \item $5H - 2E_i - 2E_j - 2E_k - 2E_l - 2E_m - E_n,$ $5H - 3E_i - 2E_j - 2E_k - E_l - E_m - E_n - E_p$.
            \item $6H - 3E_i - 3E_j - 2E_k - 2E_l - 2E_m - E_n - E_p$.
            \item $7H - 3E_i - 3E_j - 3E_k - 3E_l - 2E_m - 2E_n - E_p,$ $7H - 4E_i - 3E_j - 2E_k - 2E_l - 2E_m - 2E_n - 2E_p$.
            \item $8H - 4E_i - 3E_j - 3E_k - 3E_l - 3E_m - 2E_n - 2E_p$.
            \item $9H - 4E_i - 4E_j - 3E_k - 3E_l - 3E_m - 3E_n - 3E_p$, where $1 \leq i, j, k, l, m, n, p \leq 7$ are distinct indices.
        \end{enumerate}
        \item If $C^2 = 6$, then $C$ is big and nef. Using Kawamata-Viehweg vanishing and Riemann-Roch we can conclude that $|C+K_{\underline{X}}|$ defines a pencil. Also we have $$-K_{\underline{X}} \cdot (C+K_{\underline{X}}) = 2, \quad (C + K_{\underline{X}})^2 = 0$$ Hence $|C + K_{\underline{X}}|$ defines a birational morphism to $\PP^1$ such that $C + K_{\underline{X}}$ linearly equivalent to a smooth conic $F$, so we conclude that $C \sim -K_{\underline{X}} + F$. In particular in this case, $C$ is ample. For different choices of conics, $C$ is one of the following forms. Using \ref{prog} we get:
        \begin{enumerate}
            \item $4H - 2E_i - E_j - E_k - E_l - E_m - E_n - E_p$.
            \item $5H - 2E_i - 2E_j - 2E_k - 2E_l - E_m - E_n - E_p$.
            \item $6H - 3E_i - 2E_j - 2E_k - 2E_l - 2E_m - 2E_n - E_p$.
            \item $7H - 3E_i - 3E_j - 3E_k - 2E_l - 2E_m - 2E_n - 2E_p$.
            \item $8H - 3E_i - 3E_j - 3E_k - 3E_l - 3E_m - 3E_n - 2E_p$, where $1 \leq i, j, k, l, m, n, p \leq 7$ are distinct indices.
        \end{enumerate}
        \item If $C^2 = 8$, then $C \in |-2K_{\underline{X}}|$.
    \end{enumerate}
\end{proof}

\subsection{Classification of irreducible boundaries.} We let $\Delta = \epsilon D$, where $D$ is an irreducible curve on $\underline{X}$ and $\epsilon = 1 - 1/m.$ We discuss whether $-(K_{\underline{X}}+\Delta)$ satisfies some positivity of divisors when $D$ is an irreducible curve of anticanonical degree at most $4$.

Using Lemma \ref{thelemma} for a del Pezzo Surface of degree $2$, we can conclude that a divisor $L$ is ample (nef) if and only if $L \cdot E > 0$ (resp. $L \cdot E \geq 0$) for any $(-1)$-curve $E$. 
We use \ref{prog} to calculate the intersection numbers.

\begin{itemize}
    \item \textit{Case of $-K_{\underline{X}} \cdot D = 1:$} We assume that $D$ is a $(-1)$-curve. Since, $D \cdot E = 0, 1, 2$ where $E$ is any $(-1)$-curve. We conclude that $-(K_{\underline{X}} + \epsilon D)$ is big and nef for $\epsilon = \frac{1}{2}.$

    \item \textit{Case of $-K_{\underline{X}} \cdot D = 2:$} We assume that $D$ is a smooth fiber of a conic fibration. Then $D \cdot E = 0, 1, 2$, where $E$ is any $(-1)$-curve. So, $-(K_{\underline{X}} + \epsilon D)$ is nef for $\epsilon = \frac{1}{2},$ but not big.

    Now assume that $D$ is member of $|-K_{\underline{X}}|$. Since $D\cdot E = 1,$ we conclude that $-(K_{\underline{X}} + \epsilon D)$ is ample for all $m \in \mathbb{Z}_{\geq 1}$.
    \item \textit{Case of $-K_{\underline{X}} \cdot D = 3:$} We assume that $D$ is pullback of a hyperplane class in $\PP^2$. Then $D \cdot E = 0,1,2$ or $3$ for any $(-1)$-curve $E$. Hence we conclude that $-(K_{\underline{X}} + \epsilon D)$ is not nef for any $\epsilon > 0.$

    Now assume that $D$ is pullback of the anticanonical class of a smooth cubic surface. In this case $D \cdot E = 0, 1, 2, 3$ for any $(-1)$-curve $E$. Hence we conclude that $-(K_{\underline{X}} + \epsilon D)$ is not nef for any $\epsilon > 0$.
    \item \textit{Case of $-K_{\underline{X}} \cdot D = 4:$} We assume that $D^2 = 2$. In this case $D$ is pullback of a divisor in the $(1,1)$-class of $\PP^1 \times \PP^1$. In this case $D \cdot E = 0, 1, 2, 3, 4$, where $E$ is any $(-1)$-curve. Hence we conclude that $-(K_{\underline{X}} + \epsilon D)$ is not nef for any $\epsilon > 0.$

    Next we assume that $D^2 = 4.$ In this case $D$ is pullback of a member in the anticanonical class of a smooth quartic. Then $D \cdot E = 0, 1, 2, 3$ or $4$, where $E$ is any $(-1)$-curve. Hence, $-(K_{\underline{X}} + \epsilon D)$ is not nef for any $\epsilon > 0$.

    Next we assume that $D^2 = 6.$ In this case $D+K_{\underline{X}}$ is linearly equivalent to a smooth conic. Then $D \cdot E = 1, 2$ or $3$, where $E$ is any $(-1)$-curve. Hence we conclude that $-(K_{\underline{X}} + \epsilon D)$ is not nef for any $\epsilon > 0$.

    We assume that $D \in |-2K_{\underline{X}}|.$ Since $D \cdot E = 2,$ for any $(-1)$-curve $E$, we conclude that $-(K_{\underline{X}} + \epsilon D)$ is nef for $\epsilon = \frac{1}{2}$ but it is not big.
\end{itemize}

In summary, we conclude:

\begin{proposition}
    Let $(X, \epsilon D)$ be a del Pezzo orbifold whose underlying surface is $\underline{X}$, a del Pezzo surface of degree $2$. Let $D$ be an irreducible curve of degree less than or equal to $4$ on $\underline{X}$ and $\epsilon = 1 - 1/m, m \in \Z_{\geq 2}$.Then $-(K_{\underline{X}} + \epsilon D)$ is ample if and only if $D$ is a member of the anticanonical class $|-K_{\underline{X}}|$, for all $0 < \epsilon < 1.$
    Also, $-(K_{\underline{X}} + \epsilon D)$ is nef for $\epsilon = \frac{1}{2}$ when $D$ is either of the following.
    \begin{itemize}
        \item a $(-1)$-curve. $-(K_{\underline{X}} + \epsilon D)$ is big in this case.
        \item is a smooth fiber of a conic fibration, it is nef but not big.
        \item is a member of $|-2K_{\underline{X}}|$, it is nef but not big.
    \end{itemize}
\end{proposition}

\section{Surfaces of degree 3: Cubic surfaces.}
\label{dp3}

A del Pezzo surface of degree $3$ is a smooth cubic surface in $\PP^3$. Let $\underline{X}$ be such a del Pezzo surface. A cubic surface is the blow up $\beta: \underline{X} \to \PP^2$ at six general points. In this case the canonical divisor is given by $$K_{\underline{X}} = -3H + E_1 + E_2 + E_3 + E_4 + E_5 + E_6,$$
where $E_i$'s $(1 \leq i \leq 6)$ are the exceptional curves under $\beta$ and $\beta^*K_{\PP^2} = -3H$, where $H$ is the hyperplane class in $\underline{X}$.
We classify possible del Pezzo orbifolds $(\underline{X}, \Delta_{\epsilon})$ where $\underline{X}$ is as above and $\Delta_{\epsilon}$ is a $\Q$-divisor. In this case, the anticanonical divisor $-K_{\underline{X}}$ has degree $3$, $\Delta_{\epsilon}$ can have at most $6$ components and each component has anticanonical degree at most $6$.

\subsection{Low degree curves on $\underline{X}$.} We first classify irreducible curves of the anticanonical degree at most $3$.

\begin{itemize}
    \item An irreducible curve $C$ of anticanonical degree $1$ is a $(-1)$-curve such that $C^2 = -1$. 
    They are one of the following forms. Using \ref{prog} we get:
    \begin{enumerate}
        \item $E_i$, where $1 \leq i \leq 6$ are the $6$ exceptional curves under $\beta$.
        \item $H - E_i - E_i$, where $1 \leq i, j \leq 6, i \neq j$.
        There are $15$ such curves.
        \item $2H - E_i - E_j - E_k - E_l - E_m$, where $1 \leq i, j, k, l, m \leq 6$ are distinct indices. 
        There are $6$ such curves.
    \end{enumerate}
    \item An irreducible curve $C$ of anticanonical degree $2$ on $\underline{X}$ is a smooth fiber of a conic fibration on $\underline{X}$. Such a $C$ satisfies $C^2 = 0$. 
    They are one of the following forms. Using \ref{prog} we get:
    \begin{enumerate}
        \item $H - E_i$, where $1 \leq i \leq 6$.
        \item $2H - E_i - E_j - E_k - E_l$, where $1 \leq i, j, k, l \leq 6$ are distinct indices.
        \item $3H - 2E_i - E_j - E_k - E_l - E_m - E_n$, where $1 \leq i,j,k,l,m,n \leq 6$ are distinct indices.
    \end{enumerate}
    \item An irreducible curve $C$ of anticanonical degree $3$ on $\underline{X}$ satisfies either $C^2 = 1$ or $C^2 = 3$ by Hodge Index theorem and adjunction formula.
    In the former case, similarly as in proposition \ref{dp2d3}, Case 1, we can conclude that $C$ is a member of the linear system of the pullback of a hyperplane class via a birational morphism to $\PP^2$. 
    They are one of the following forms. Using \ref{prog} we get:
    \begin{enumerate}
        \item $H$.
        \item $2H - E_i - E_j - E_k$, where $1 \leq i, j, k \leq 6$ are distinct indices.
        \item $3H - 2E_i - E_j - E_k - E_l - E_m$, where $1 \leq i,j,k,l,m \leq 6$ are distinct indices.
        \item $4H - 2E_i - 2E_j - 2E_k - E_l - E_m - E_n$.
        \item $5H - 2E_i - 2E_j - 2E_k - 2E_l - 2E_m - 2E_n$, where $H$ is the pullback of the hyperplane class in $\PP^2$ by $\beta$, and $1 \leq i,j,k,l,m,n \leq 6$ are distinct indices.
    \end{enumerate}
    In the later case $C$ is a member of $|-K_{\underline{X}}|$. 
\end{itemize}

Next we assume that $C$ is an irreducible curve of anticanonical degree $4$. Then by Hodge Index theorem we must have $0 \leq C^2 \leq 5$. Also we must have $2p_a(C) - 2 = -4 +C^2$, so we conclude that $C^2 = 2$ or $4$.

\begin{proposition}
    Let $\underline{X}$ be a cubic surface and $C$ is an irreducible curve of anticanonical degree $4$, then either $C^2 = 2$ or $C^2 = 4$.
    When $C^2 = 2$, there exists a birational map $\pi: \underline{X} \to \PP^1 \times \PP^1$ such that $C$ is pullback of a member in the $(1,1)$-class in $\PP^1 \times \PP^1$.
    When $C^2 = 4$, there exists a birational morphism $\pi: \underline{X} \to \underline{X}'$ to a smooth quartic del Pezzo surface such that $C$ is a member of $|-\pi^*K_{\underline{X}'}|$.
\end{proposition}

\begin{proof}
    \begin{enumerate}
        \item \textit{Case 1: $C^2 = 2.$} Similarly as in proposition \ref{dp2d4}, Case $1$, we can conclude that $C$ is pullback of a member in the $(1,1)$-class in $\PP^1 \times \PP^1$ via a birational morphism.  
        They are of the following forms. Using \ref{prog} we get:
        \begin{enumerate}
            \item $2H - E_i - E_j$, where $1 \leq i, j \leq 6, i \neq j$.
            \item $3H - 2E_i - E_j - E_k - E_l$, where $1 \leq i,j,k,l \leq 6$ are distinct indices.
            \item $4H - 2E_i - 2E_j - 2E_k - E_l - E_m$; $4H - 3E_i - E_j - E_k - E_l - E_m - E_n$.
            \item $5H - 3E_i - 2E_j - 2E_k - 2E_l - E_m - E_n$.
            \item $6H - 3E_i - 3E_j - 2E_k - 2E_l - 2E_m - 2E_n$, where $1 \leq i,j,k,l,m,n \leq 6$ are distinct indices.
        \end{enumerate}
        \item \textit{Case 2: $C^2 = 4.$} Similarly as in proposition \ref{dp2d4}, case 2, we can conclude that $C$ is must be a member of the pullback of anticanonical linear system on a smooth quartic del Pezzo surface via a birational morphism.
        They are one of the following forms. Using \ref{prog} we get:
        \begin{enumerate}
            \item $3H - E_i - E_j - E_k - E_l - E_m$, where $1 \leq i,j,k,l,m \leq 6$ are distinct indices.
            \item $4H - 2E_i - 2E_j - E_k - E_l - E_m - E_n$.
            \item $5H - 2E_i - 2E_j - 2E_k - 2E_l - 2E_m - E_n$, where $1 \leq i,j,k,l,m,n \leq 6$ are distinct indices. 
        \end{enumerate}
    \end{enumerate}
\end{proof}

Next we assume that $C$ is an irreducible divisor of anticanonical degree $5$. Then by Hodge Index theorem and adjunction formula we must have $C^2 = 3, 5$ or $7$.

\begin{proposition}
    \label{dp3d5}
    Let $\underline{X}$ be a cubic surface and $C$ is an irreducible curve of anticanonical degree $5$, then either $C^2 = 3, 5$ or $7$.
    When $C^2 = 3$, there exists a birational morphism $\pi: \underline{X} \to \underline{X}'$ to blow up of $\PP^2$ at one point such that $C$ is a member of $|2\pi^*H' - \pi^* E|$ where $H'$ is the pullback of the hyperplane class in $\PP^2$ via a blow up $\beta': \underline{X}' \to \PP^2$ at one point, and $E$ is the unique exceptional divisor of $\beta'$.
    When $C^2 = 5$, there exists a birational morphism $\pi: \underline{X} \to \underline{X}'$ to a smooth quintic del Pezzo surface such that $C$ is a member of $|-\pi^*K_{\underline{X}'}|$.
    When $C^2 = 7$, then $C$ is ample and one can find a smooth conic $F$ such that $C \sim -K_{\underline{X}} + F$.
\end{proposition}

\begin{proof}
    \begin{enumerate}
        \item \textit{Case 1: $C^2 = 3.$} In this case $2C$ is big and nef. Hence using Kawamata-Viehweg vanishing and Riemann-Roch we can conclude that $|2C+K_{\underline{X}}|$ defines a pencil. We also have $$C \cdot (2C + K_{\underline{X}}) = 1, \quad  -K_{\underline{X}} \cdot (2C + K_{\underline{X}}) = 7, \quad (2C + K_{\underline{X}})^2 = -5.$$
        Then by Zariski decomposition we can write $2C + K_{\underline{X}} \sim F + \sum_i a_i N_i \geq 0$ where $a_i > 0$. Since $(2C + K_{\underline{X}}) \cdot N_i < 0$, we can conclude that $N_i$'s are disjoint $(-1)$-curves contracted by $C$. Then we must have 
        $$-K_{\underline{X}} \cdot F + \sum_i a_i = 7, \quad F \cdot N_i = a_i - 1, \quad F^2 - \sum_i a_i^2 + \sum_i 2a_i F \cdot N_i = -5.$$
        Thus we conclude $$F^2 + \sum_i a_i^2 - \sum_i 2a_i = -5$$ By Zariski decomposition $F \cdot N_i = 0$ for all $i$, hence we have $-K_{\underline{X}} \cdot F = 2, F^2 = 0$ and $a_i = 1$ for all $i$. 
        Hence the number of $N_i$'s is $5$. Let $\pi: \underline{X} \to \underline{X}'$ be a birational contraction of $N_i$'s and denote the pushforward of $C$ and $F$ by $C'$ and $F'$. Then on $\underline{X}'$ we have $$2C' + K_{\underline{X}'} \sim F'$$ 
        and we also have $-K_{\underline{X}'} \cdot C' = 5$. This is only possible when $\underline{X}'$ is blow up of $\PP^2$ at one point. Let $H'$ be the pullback of the hyperplane class in $\PP^2$ via a blow up $\beta': \underline{X}' \to \PP^2$ at one point, and $E$ be the unique exceptional divisor of $\beta'$. Then $C' \sim 2H' - E$ and $C$ is a member of $|2\pi^*H' - \pi^*E|$.
        They are one of the following forms. Using \ref{prog} we get:
        \begin{enumerate}
            \item $2H - E_i$, where $1 \leq i \leq 6$.
            \item $3H - 2E_i - E_j - E_k$, where $1 \leq i,j,k \leq 6$ are distinct indices.
            \item $4H - 2E_i - 2E_j - 2E_k - E_l$; $4H - 3E_i - E_j - E_k - E_l - E_m$.
            \item $5H - 3E_i - 2E_j - 2E_k - 2E_l - E_m$; where $1 \leq i,j,k,l,m,n \leq 6$ are distinct indices.
            \item $6H - 3E_i - 3E_j - 3E_k - 2E_l - E_m - E_n$; $6H - 4E_i - 2E_j - 2E_k - 2E_l - 2E_m - E_n$.
            \item $7H - 4E_i - 3E_j - 3E_k - 2E_l - 2E_m - 2E_n$.
            \item $8H - 4E_i - 3E_j - 3E_k - 3E_l - 3E_m - 3E_n$, where $1 \leq i,j,k,l,m,n \leq 6$ are distinct indices.
        \end{enumerate}
        \item \textit{Case 2: $C^2 = 5.$} In this case $C+K_{\underline{X}}$ is linearly equivalent to the disjoint union of two $(-1)$-curves. Let $\pi: \underline{X} \to \underline{X}'$ be the birational contraction of these $(-1)$-curves. Then $\underline{X}'$ is a degree $5$ del Pezzo surface and $C$ is a member of $|-\pi^*K_{\underline{X}'}|$. 
        They are one of the following forms. Using \ref{prog} we get:
        \begin{enumerate}
            \item $3H - E_i - E_j - E_k - E_l$, where $1 \leq i,j,k,l \leq 6$ are distinct indices.
            \item $4H - 2E_i - 2E_j - E_k - E_l - E_m$, where $1 \leq i,j,k,l,m \leq 6$ are distinct indices.
            \item $5H - 2E_i - 2E_j - 2E_k - 2E_l - 2E_m$; $5H - 3E_i - 2E_j - 2E_k - E_l - E_m - E_n$.
            \item $6H - 3E_i - 3E_j - 2E_k - 2E_l - 2E_m - E_n$.
            \item $7H - 3E_i - 3E_j - 3E_k - 3E_l - 2E_m - 2E_n$, where $1 \leq i,j,k,l,m,n \leq 6$ are distinct indices. 
        \end{enumerate}
        \item \textit{Case 3: $C^2 = 7.$} Similarly as in proposition \ref{dp2d4}, Case 3, we can conclude that $C + K_{\underline{X}}$ is linearly equivalent to a smooth fiber $F$ of a conic fibration $f: \underline{X} \to \PP^1$. So we conclude that $C \sim -K_{\underline{X}} + F$. In particular, in this case $C$ is ample.
        They are one of the following forms.  Using \ref{prog} we get:
        \begin{enumerate}
            \item $4H - 2E_i - E_j - E_k - E_l - E_m - E_n$.
            \item $5H - 2E_i - 2E_j - 2E_k - 2E_l - E_m - E_n$.
            \item $6H - 3E_i - 2E_j - 2E_k - 2E_l - 2E_m - 2E_n$, where $1 \leq i,j,k,l,m,n \leq 6$ are distinct indices.
        \end{enumerate}
    \end{enumerate}
\end{proof}

\begin{proposition}
    \label{dp3d6}
    Let $\underline{X}$ be a cubic surface and $C$ is an irreducible curve of anticanonical degree $6$, then either $C^2 = 4, 6, 8, 10$ or $12$.
    When $C^2 = 4$, there exists a birational morphism $\pi: \underline{X} \to \underline{X}'$ to $\PP^1 \times \PP^1$ such that $C$ is pullback of a member of $(2,1)$ or $(1,2)$-class in $\PP^1 \times \PP^1$. 
    When $C^2 = 6$, there exists a birational morphism $\pi: \underline{X} \to \underline{X}'$ to a smooth sextic del Pezzo surface such that $C$ is a member of $|-\pi^*K_{\underline{X}'}|$.
    When $C^2 = 8$, there exists a birational morphism $\pi: \underline{X} \to \underline{X}'$ to a del Pezzo surface of degree $4$ such that $C$ is a member of $|-\pi^*K_{\underline{X}'} + \pi^*F|$, for some smooth conic $F$. 
    When $C^2 = 10$, there exists a birational morphism $\pi: \underline{X} \to \underline{X}'$ to $\PP^2$ such that $C$ is a member of $|-K_{\underline{X}} + \pi^* H'|$, where $H'$ is the hyperplane class in $\PP^2$.
    When $C^2 = 12$, then $C$ is a member of $|-2K_{\underline{X}}|$.
\end{proposition}

\begin{proof}
    \begin{enumerate}
        \item \textit{Case 1: $C^2 = 4.$} In this case $2C$ is big and nef. Using Kawamata-Viehweg vanishing and Riemann-Roch we can conclude that $2C + K_{\underline{X}}$ is effective. We also have $$C \cdot (2C + K_{\underline{X}}) = 2, \quad -K_{\underline{X}} \cdot (2C + K_{\underline{X}}) = 9, \quad (2C + K_{\underline{X}})^2 = -5.$$
        After writing the Zariski decomposition of the divisor $2C + K_{\underline{X}}$, by similar calculations as in proposition \ref{dp3d5}, Case 1, we can conclude that the number of $N_i$'s is $5$. Let $\pi: \underline{X} \to \underline{X}'$ be a birational contraction of $N_i$'s.
        This is only possible when $\underline{X}'$ is blow-up of $\PP^2$ at one point or is $\PP^1 \times \PP^1$. Let $H'$ be pullback of a hyperplane class in $\PP^2$ via a blow up $\beta': \underline{X}' \to \PP^2$ and $E$ be the unique exceptional divisor of $\beta'$, then $2C' \sim 5H' - 3E$.  But since $C'$ is irreducible and the divisor $5H - 3E$ is not divisible by $2$, this is not possible.
        Hence $\underline{X}'$ must be $\PP^1 \times \PP^1$, and $C'$ must be a member of a $(2,1)$ or $(1,2)$-class. Hence $C$ is pullback of a member either in $(2,1)$-class or in $(1,2)$-class in $\PP^1 \times \PP^1$ via a birational morphism. 
        For different contraction maps, they are one of the following forms. Using \ref{prog} we get:
        \begin{enumerate}
            \item $2H$.
            \item $3H - 2E_i - E_j$, where $1 \leq i,j \leq 6, i \neq j$.
            \item $4H - 2E_i - 2E_j - 2E_k$; $4H - 3E_i - E_j - E_k - E_l$, where $1 \leq i,j,k,l \leq 6$ are distinct indices.
            \item $5H - 3E_i - 2E_j - 2E_k - 2E_l$; $5H - 4E_i - E_j - E_k - E_l - E_m - E_n$.
            \item $6H - 3E_i - 3E_j - 3E_k - 2E_l - E_m$; $6H - 4E_i - 2E_j - 2E_k - 2E_l - 2E_m$.
            \item $7H - 4E_i - 3E_j -3E_k - 3E_l - E_m - E_n$; $7H - 5E_i - 2E_j - 2E_k - 2E_l - 2E_m - 2E_n$.
            \item $8H - 4E_i - 4E_j - 4E_k - 2E_l - 2E_m - 2E_n$; $8H - 5E_i - 3E_j - 3E_k - 3E_l - 2E_m - 2E_n$.
            \item $9H - 5E_i - 4E_j - 3E_k - 3E_l - 3E_m - 3E_n$, where $1 \leq i,j,k,l,m,n \leq 6$ are distinct indices.
            \item $10H - 4E_1 - 4E_2 - 4E_3 - 4E_4 - 4E_5 - 4E_6$. 
        \end{enumerate}
        \item \textit{Case 2: $C^2 = 6.$} In this case $C + K_{\underline{X}}$ is linearly equivalent to disjoint union of three $(-1)$-curves. Let $\pi: \underline{X} \to \underline{X}'$ be the contraction of these curves to a smooth sextic del Pezzo surface. Then $C$ is pullback of a member in the anticanonical linear system of the smooth sextic del Pezzo surface.
        They are one of the following forms. Using \ref{prog} we get:
        \begin{enumerate}
            \item $3H - E_i - E_j - E_k$, where $1 \leq i,j,k \leq 6$ are distinct indices.
            \item $4H - 2E_i - 2E_j - E_k - E_l$, where $1 \leq i,j,k,l \leq 6$ are distinct indices.
            \item $5H - 3E_i - 2E_j - 2E_k - E_l - E_m$, where $1 \leq i,j,k,l,m \leq 6$ are distinct indices.
            \item $6H - 3E_i - 3E_j - 2E_k - 2E_l - 2E_m$; $6H - 3E_i - 3E_j - 3E_k - E_l - E_m - E_n$; $6H - 4E_i - 2E_j - 2E_k - 2E_l - E_m - E_n$.
            \item $7H - 4E_i - 3E_j - 3E_k - 2E_l - 2E_m - E_n$.
            \item $8H - 4E_i - 4E_j - 3E_k - 3E_l - 2E_m - 2E_n$.
            \item $9H - 4E_i - 4E_j - 4E_k - 3E_l - 3E_m - 3E_n$, where $1 \leq i,j,k,l,m,n \leq 6$ are distinct indices.
        \end{enumerate}
        \item \textit{Case 3: $C^2 = 8.$} In this case $C$ is big and nef. Using Kawamata-Viehweg vanishing and Riemann-Roch we can conclude that $C + K_{\underline{X}}$ is effective. We also have $$(C + K_{\underline{X}})^2 = -1, \quad -K_{\underline{X}} \cdot (C + K_{\underline{X}}) = 3, \quad C \cdot (C + K_{\underline{X}}) = 2.$$
        After writing the Zariski decomposition of the divisor $C + K_{\underline{X}}$, using similar calculations as in proposition \ref{dp3d5}, Case 1, we can conclude that the number of $N_i$'s is one. Let $\pi: \underline{X} \to \underline{X}'$ be a birational contraction of this $N_i$.
        This is only possible when $\underline{X}'$ is a del Pezzo surface of degree $4$. Let $C'$ denote the pushforward of $C$ in $\underline{X}'$. Then $C'$ is a curve with anticanonical degree $6$ and $C'^2 = 8$ on a smooth del Pezzo surface of degree $4$. As described below in proposition \ref{dp4d6}, case 3, $C'$ is linearly equivalent to $-K_{\underline{X}'} + F$ for some smooth conic $F$. Therefore, $C$ is a member of $|-\pi^*K_{\underline{X}'} + \pi^*F|$.
        They are one of the following forms. Using \ref{prog} we get:
        \begin{enumerate}
            \item $4H - 2E_i - E_j - E_k - E_l - E_m$, where $1 \leq i,j,k,l,m \leq 6$ are distinct indices.
            \item $5H - 2E_i - 2E_j - 2E_k - 2E_l - E_m$; $5H - 3E_i - 2E_j - E_k - E_l - E_m - E_n$.
            \item $6H - 3E_i - 3E_j - 2E_k - 2E_l - E_m - E_n$.
            \item $7H - 3E_i - 3E_j - 3E_k - 3E_l - 2E_m - E_n$; $7H - 4E_i - 3E_j - 2E_k - 2E_l - 2E_m - 2E_n$.
            \item $8H - 4E_i - 3E_j - 3E_k - 3E_l - 3E_m - 2E_n$, where $1 \leq i,j,k,l,m,n \leq 6$ are distinct indices.
        \end{enumerate}
        \item \textit{Case 4: $C^2 = 10.$} Using Kawamata-Viehweg vanishing and Riemann-Roch $C + K_{\underline{X}}$ must be effective. We also have $-K_{\underline{X}} \cdot (C + K_{\underline{X}}) = 3$ and $(C + K_{\underline{X}})^2 = 1$. Hence $C + K_{\underline{X}}$ is a member of the linear system of the pullback of a hyperplane class via a birational morphism $\pi: \underline{X} \to \underline{X}'$ to $\PP^2$. Hence $C \in |-K_{\underline{X}} + \pi^*H'|$, where $H'$ is the hyperplane class in $\PP^2$.
        They are one of the following forms. Using \ref{prog} we get:
        \begin{enumerate}
            \item $4H - E_i - E_j - E_k - E_l - E_m - E_n$.
            \item $5H - 2E_i - 2E_j - 2E_k - E_l - E_m - E_n$.
            \item $6H - 3E_i - 2E_j - 2E_k - 2E_l - 2E_m - E_n$.
            \item $7H - 3E_i - 3E_j - 3E_k - 2E_l - 2E_m - 2E_n$.
            \item $8H - 3E_i - 3E_j - 3E_k - 3E_l - 3E_m - 3E_n$, where $1 \leq i,j,k,l,m,n \leq 6$ are distinct indices.
        \end{enumerate}
        \item \textit{Case 5: $C^2 = 12.$} In this case $C \in |-2K_{\underline{X}}|$.
    \end{enumerate}
\end{proof}

\subsection{Classification of irreducible boundaries.} We let $\Delta = \epsilon D$, where $D$ is an irreducible curve on $\underline{X}$ and $\epsilon = 1 - 1/m.$ We discuss whether $-(K_{\underline{X}}+\Delta)$ satisfies some positivity of divisors when $D$ is an irreducible curve of anticanonical degree at most $6$.

Using Lemma \ref{thelemma} for a del Pezzo Surface of degree $3$, we can conclude that a divisor $L$ is ample (nef) if and only if $L \cdot E > 0$ (resp. $L \cdot E \geq 0$) for any $(-1)$-curve $E$. 
We use \ref{prog} to calculate the intersection numbers.

\begin{itemize}
    \item \textit{Case of $-K_{\underline{X}} \cdot D = 1$:} We assume that $D$ is a $(-1)$-curve. Since $D \cdot E = 0$ or $1$ for any $(-1)$-curve $E$, we conclude that $-(K_{\underline{X}} + \epsilon D)$ is ample for any $\epsilon = 1 - 1/m$ where $m \in \Z_{\geq 1}$. 
    In this case, $-K_{\underline{X}} + D$ is big and nef but not ample. 

    \item \textit{Case of $-K_{\underline{X}} \cdot D = 2$:} We assume that $D$ is a smooth conic. In this case $D \cdot E = 0, 1$ or $2$ for any $(-1)$-curve $E$, so we conclude that $-(K_{\underline{X}} + \epsilon D)$ is big and nef for $\epsilon = \frac{1}{2}$.

    \item \textit{Case of $-K_{\underline{X}} \cdot D = 3$:} First we assume that $D^2 = 1$. In this case $D$ is pullback of a hyperplane class on $\PP^2$. Since pushforward of a $(-1)$-curve on $\PP^2$ is either $0$, a line, or a conic, we conclude that $D \cdot E = 0, 1$ or $2$ for any $(-1)$-curve $E$. 
    This imply that $-(K_{\underline{X}} + \epsilon D)$ is big and nef for $\epsilon = \frac{1}{2}$.

    Next assume that $D$ is linearly equivalent to $-K_{\underline{X}}$. In this case $-(K_{\underline{X}} + \epsilon D)$ is ample for any $\epsilon = 1 - 1/m$ where $m \in \Z_{\geq 1}$ and $-(K_{\underline{X}}+D)$ is trivial.

    \item \textit{Case of $-K_{\underline{X}} \cdot D = 4$:} First we assume that $D^2 = 2$ so that $D$ is the pullback of the hyperplane class of a smooth quadric. The pushforward of a $(-1)$-curve to the quadric is either contracted or has class $(1,0), (1,1)$ or their involutes. So $D \cdot E = 0, 1, 2$ or $3$. This imply that $-(K_{\underline{X}} + \epsilon D)$ is not nef for any $\epsilon$.

    Next we assume that $D^2 = 4$ so that $D$ is pullback of the anticanonical divisor on a quartic del Pezzo surface. Then $D \cdot E = 0, 1$ or $2$. This imply that $-(K_{\underline{X}} + \epsilon D)$ is nef but not big for $\epsilon = \frac{1}{2}$.

    \item \textit{Case of $-K_{\underline{X}} \cdot D = 5$:} First we assume that $D^2 = 3.$ Then $D$ is the pullback of a divisor of the form $2H' - E$ as explained above, where $H'$ is the pullback of the hyperplane class from $\PP^2$ via a blow up $\beta':\underline{X}' \to \PP^2$ at one point and $E$ is the exceptional curve of $\beta'$. In this case $D \cdot F = 0, 1, 2, 3$ or $4$ for any $(-1)$-curve $F$. Hence $-(K_{\underline{X}} + \epsilon D)$ is not nef for any $\epsilon$.
    
    Next we assume that $D^2 = 5.$ Then $D$ is pullback of a member of the anticanonical class of a del Pezzo surface of degree $5$. In this case $D \cdot E = 0, 1, 2$ or $3$ for any $(-1)$-curve $E$. Hence $-(K_{\underline{X}} + \epsilon D)$ is not nef for any $\epsilon$.
    
    Next we assign that $D^2 = 7$. Then $D$ is linearly equivalent to $-K_{\underline{X}} + F$ where $F$ is a smooth conic. In this case $D \cdot E = 1, 2$ or $3$ for any $(-1)$-curve $E$. Hence $-(K_{\underline{X}} + \epsilon D)$ is not nef for any $\epsilon$. 
    \item \textit{Case of $-K_{\underline{X}} \cdot D = 6$:} First we assume that $D^2 = 4$. In this case $D \cdot E$ can be as large as $5$ for some $(-1)$-curve $E$. Hence $-(K_{\underline{X}} + \epsilon D)$ can not be nef for any $\epsilon$.

    Next we assume that $D^2 = 6$ so that $D$ is pullback of a member in the anticanonical class of a smooth sextic del Pezzo surface. In this case $D \cdot E = 0, 1, 2, 3$ or $4$ for any $(-1)$-curve $E$. Hence $-(K_{\underline{X}} + \epsilon D)$ is not nef for any $\epsilon$.

    Next we assume that $D^2 = 8$. In this case $D \cdot E$ can be $0, 1, 2, 3$ or $4$ for any $(-1)$-curve $E$. Hence $-(K_{\underline{X}} + \epsilon D)$ is not nef for any $\epsilon$.

    Next we assume that $D^2 = 10$. In this case $D \cdot E = 1, 2$ or $3$ for any $(-1)$-curve $E$. Hence we conclude that $-(K_{\underline{X}} + \epsilon D)$ is not nef for any $\epsilon$.

    Next we assume that $D^2 = 12$ so that $D \in |-2K_{\underline{X}}|$. Then we have $D \cdot E = 2$ for all $(-1)$-curve $E$. Hence $-(K_{\underline{X}} + \epsilon D)$ is nef for $\epsilon = \frac{1}{2}$, but not big.
\end{itemize}

In summary, we conclude:

\begin{proposition}
    Let $(X, \epsilon D)$ be a del Pezzo orbifold whose underlying surface $\underline{X}$ is a smooth cubic surface. Let $D$ be an irreducible curve of degree less than or equal to $6$, and $\epsilon = 1 - 1/m, m \in \Z_{\geq 2}$. Then $-(K_{\underline{X}} + \epsilon D)$ is ample if and only if
    \begin{itemize}
        \item $D$ is a $(-1)$-curve and $\epsilon > 0$.
        \item $D$ is a member of the anticanonical class and $\epsilon > 0$.
    \end{itemize} 
    Also $-(K_{\underline{X}} + \Delta)$ is nef when
    \begin{itemize}
        \item $D$ is pullback of a hyperplane class in $\PP^2$ and $\epsilon = \frac{1}{2}$, in this case $-(K_{\underline{X}} + \Delta)$ is also big.
        \item $D$ is a smooth fiber of a conic fibration and $\epsilon = \frac{1}{2}$, in this case $-(K_{\underline{X}} + \Delta)$ is also big.
        \item $D$ is pullback of a anticanonical divisor on a quartic del Pezzo surface and $\epsilon = \frac{1}{2}$.
        \item $D$ is a member of $|-2K_{\underline{X}}|$ and $\epsilon = \frac{1}{2}$. 
    \end{itemize}
\end{proposition}

\section{Surfaces of degree 4.}
\label{dp4}

A del Pezzo surface of degree $4$ is intersection of two quadrics in $\PP^4.$ Let $\underline{X}$ be such a del Pezzo surface. $\underline{X}$ is the blow up $\beta: \underline{X} \to \PP^2$ at $5$ general points. In this case the canonical divisor is given by 
$$K_{\underline{X}} = -3H + E_1 + E_2 + E_3 + E_4+ E_5,$$ where $E_i$'s $(1 \leq i \leq 5)$ are the exceptional curves of $\beta$ and $\beta^*K_{\PP^2} = -3H$, where $H$ is the hyperplane class in $\underline{X}$. 
We classify possible del Pezzo orbifolds $(\underline{X}, \Delta_{\epsilon})$ where $\underline{X}$ is as above and $\Delta_{\epsilon}$ is a $\Q$-divisor. In this case, the anticanonical divisor $-K_{\underline{X}}$ has degree $4$, $\Delta_{\epsilon}$ can have at most $8$ components and each component has anticanonical degree at most $8$.

\subsection{Low degree curves on $\underline{X}$.} We classify all irreducible curves with anticanonical degree less than $5$.

\begin{itemize}
    \item An irreducible curve $C$ of anticanonical degree $1$ is a $(-1)$-curve such that $C^2 = -1.$ There are $16$ such curves, which are of the following form. Using \ref{prog} we get:
    \begin{enumerate}
        \item $E_i$ for $i = 1, \ldots, 5$.
        \item $H - E_i - E_j$, where $1 \leq i, j \leq 5, i \neq j$. There are $10$ such curves.
        \item $2H - E_1 - E_2 - E_3 - E_4 - E_5$. There is one such curve.
    \end{enumerate}
    \item Let $C$ is be irreducible curve with anticanonical degree $2$. Then by Hodge Index theorem and adjunction formula, we must have $C^2 = 0$.
    In this case $|C|$ is a basepoint-free pencil. Hence it defines a morphism $f: \underline{X} \to \PP^1.$ Such $C$ is a smooth fiber of a conic fibration $f$ to $\PP^1$. They are one of the following forms. Using \ref{prog} we get:
    \begin{enumerate}
        \item $H - E_i$, where $1 \leq i \leq 5$.
        \item $2H - E_i - E_j - E_k - E_l$, where $1 \leq i, j, k, l \leq 5$ are distinct indices.
    \end{enumerate}
    \item Let $C$ be an irreducible curve with anticanonical degree $3$. Then using Hodge Index theorem and adjunction, $C^2 = 1.$ 
    Similarly as in proposition \ref{dp2d3}, Case 1, we can conclude that $C$ is a member of the pullback of a hyperplane class in $\PP^2$ via a birational morphism.
    They are of the following forms. Using \ref{prog} we get:
    \begin{enumerate}
        \item $H$.
        \item $2H - E_i - E_j - E_k$, where $1 \leq i, j, k \leq 5$ are distinct indices.
        \item $3H - 2E_i - E_j - E_k - E_l - E_m$, where $1 \leq i, j, k, l, m \leq 5$ are distinct indices.
    \end{enumerate}
    \item Let $C$ be an irreducible curve with anticanonical degree $4$. Then using Hodge Index theorem and adjunction formula, we must have $C^2 = 2$ or $4$.
    \begin{enumerate}
        \item \textit{Case 1: $C^2 = 2.$} Similarly as in proposition \ref{dp2d4}, Case 1, we can conclude that $C$ is pullback of a member in $(1,1)$ class of $\PP^1 \times \PP^1$. They are one of the following forms. Using \ref{prog} we get:
        \begin{enumerate}
            \item $2H - E_i - E_j$, where $1 \leq i, j \leq 5, i \neq j$.
            \item $3H - 2E_i - E_j - E_k - E_l$, where $1 \leq i, j, k, l \leq 5$ are distinct indices.
            \item $4H - 2E_i - 2E_j - 2E_k - E_l - E_m$, where $1 \leq i, j, k, l, m \leq 5$ are distinct indices.
        \end{enumerate}
        \item \textit{Case 2: $C^2 = 4.$} Then $C$ is a member of the anticanonical class $|-K_{\underline{X}}|$.
    \end{enumerate}
\end{itemize}

Now assume that $C$ is an irreducible curve with anticanonical degree $5$. Then using Hodge Index theorem and adjunction formula, we must have $C^2 = 3$ or $5$.

\begin{proposition}
    Let $\underline{X}$ be a del Pezzo surface of degree $4$ and $C$ be an irreducible curve of anticanonical degree $5$. Then we have either $C^2 = 3$ or $5$.
    When $C^2 = 3$, there exists a birational morphism $\pi: \underline{X} \to \underline{X}'$ to blow-up of $\PP^2$ at one point such that $C$ is a member of $|2\pi^*H' - \pi^*E|$ where $H'$ is the pullback of a hyperplane class in $\PP^2$ via a blow up $\beta': \underline{X}' \to \PP^2$ at one point and $E$ is the unique exceptional divisor of $\beta'$.
    When $C^2 = 5$, there exists a birational morphism $\pi: \underline{X} \to \underline{X}'$ to a smooth quintic del Pezzo surface such that $C$ is a member of $|-\pi^*K_{\underline{X}'}|$.
\end{proposition}

\begin{proof}
    \begin{enumerate}
        \item \textit{Case 1: $C^2 = 3.$} This case is similar to the proof of proposition \ref{dp3d5}. Using \ref{prog} we get:
        \begin{enumerate}
            \item $2H - E_i$, where $1 \leq i \leq 5$.
            \item $3H - 2E_i - E_j - E_k$, where $1 \leq i, j, k \leq 5$ are distinct indices.
            \item $4H - 2E_i - 2E_j - 2E_k - E_l$, where $1 \leq i, j, k, l \leq 5$ are distinct indices.
            \item $4H - 3E_i - E_j - E_k - E_l - E_m$.
            \item $5H - 3E_i - 2E_j - 2E_k - 2E_l - E_m$, where $1 \leq i, j, k, l, m \leq 5$ are distinct indices.
        \end{enumerate}
        \item \textit{Case 2: $C^2 = 5.$} This case is also similar to the proof of proposition \ref{dp3d5}. Using \ref{prog} we get:
        \begin{enumerate}
            \item $3H - E_i - E_j - E_k - E_l$.
            \item $4H - 2E_i - 2E_j - E_k - E_l - E_m,$ where $1 \leq i, j, k, l, m \leq 5$ are distinct indices.
            \item $5H - 2E_i - 2E_j - 2E_k - 2E_l$, where $1 \leq i, j, k, l \leq 5$ are distinct indices.
        \end{enumerate}
    \end{enumerate}
\end{proof}

Next assume that $C$ is an irreducible curve with anticanonical degree $6$. Then using Hodge Index theorem and adjunction formula, we must have $C^2 = 4, 6$ or $8$.

\begin{proposition}
    \label{dp4d6}
    Let $\underline{X}$ be a del Pezzo surface of degree $4$ and $C$ is an irreducible curve of anticanonical degree $6$. Then we have either $C^2=4, 6$ or $8$.
    When $C^2 = 4$, there exists a birational morphism $\pi: \underline{X} \to \underline{X}'$ to $\PP^1 \times \PP^1$ such that $C$ is pullback of a member of $(2,1)$ or $(1,2)$-class in $\PP^1 \times \PP^1$. 
    When $C^2 = 6$, there exists a birational morphism $\pi: \underline{X} \to \underline{X}'$ to a smooth sextic del Pezzo surface such that $C$ is a member of $|-\pi^*K_{\underline{X}'}|$.
    When $C^2 = 8$, then $C$ is ample and one can find a smooth conic $F$ such that $C \sim -K_{\underline{X}} + F$.
\end{proposition}

\begin{proof}
    \begin{enumerate}
        \item \textit{Case 1: $C^2 = 4.$} This case is similar to the proof of proposition \ref{dp3d6}, case 1. Using \ref{prog} we get the following forms:
        \begin{enumerate}
            \item $2H$.
            \item $4H - 2E_i - 2E_j - 2E_k$, where $1 \leq i, j, k \leq 5$ are distinct indices.
            \item $5H - 3E_i - 2E_j - 2E_k - 2E_l$, where $1 \leq i, j, k, l \leq 5$ are distinct indices.
            \item $6H - 4E_i - 2E_j - 2E_k - 2E_l - 2E_m$.
            \item $3H - 2E_i - E_j$, where $1 \leq i, j \leq 5, i \neq j$.
            \item $4H - 3E_i - E_j - E_k - E_l$, where $1 \leq i, j, k, l \leq 5$ are distinct indices.
            \item $6H - 3E_i - 3E_j - 3E_k - 2E_l - E_m$, where $1 \leq i, j, k, l, m \leq 5$ are distinct indices.
        \end{enumerate}
        \item \textit{Case 2: $C^2 = 6.$} This case is similar to the proof of proposition \ref{dp3d6}, case 2. Using \ref{prog} we get the following forms:
        \begin{enumerate}
            \item $3H - E_i - E_j - E_k$, where $1 \leq i, j, k \leq 5$ are distinct indices.
            \item $4H - 2E_i - 2E_j - E_k - E_l$, where $1 \leq i, j, k, l \leq 5$ are distinct indices.
            \item $5H - 3E_i - 2E_j - 2E_k - E_l - E_m$.
            \item $6H - 3E_i - 3E_j - 2E_k - 2E_l - 2E_m$, where $1 \leq i, j, k, l, m \leq 5$ are distinct indices.
        \end{enumerate}
        \item \textit{Case 3: $C^2 = 8.$} In this case similarly as in proposition \ref{dp3d5}, Case 3, one can prove that $C + K_{\underline{X}}$ is linearly equivalent a smooth conic $F$, so we conclude that $C \sim -K_{\underline{X}} + F$. In particular in this case, $C$ is ample. Using \ref{prog} we get two type of conics such that
        $C = 4H - 2E_i - E_j - E_k - E_l - E_m$ or, $5H - 2E_i - 2E_j - 2E_k - 2E_l - E_m$, where $1 \leq i, j, k, l, m \leq 5$ are distinct indices.
    \end{enumerate}
\end{proof}

Next let us assume $C$ is an irreducible curve with anticanonical degree $7$. Then using Hodge Index theorem and adjunction, $C^2 = 5,7,9,11.$

\begin{proposition}
    \label{dp4d7}
    Let $\underline{X}$ be a del Pezzo surface of degree $4$ and $C$ be an irreducible curve of anticanonical degree $7$. Then we have either $C^2 = 5, 7, 9$ or, $11$.
    When $C^2 = 5$, there exists a birational morphism $\pi: \underline{X} \to \underline{X}'$ to the blow up of $\PP^2$ at a point such that $C$ is a member of $|3\pi^*H' - 2\pi^*E|$ where $H'$ is the pullback of the hyperplane class from $\PP^2$ via a blow up $\beta':\underline{X}' \to \PP^2$ at one point and $E$ is the unique exceptional divisor of $\beta'$.
    When $C^2 = 7$, there exists a birational morphism $\pi: \underline{X} \to \underline{X}'$ to a smooth del Pezzo surface of degree $7$ such that $C$ is a member of $|-\pi^*K_{\underline{X}'}|.$
    When $C^2 = 9$, there exists a birational morphism $\pi: \underline{X} \to \underline{X}'$ to a smooth del Pezzo surface of degree $5$ such that $C$ is a member of $|-\pi^*K_{\underline{X}'} + \pi^*F|$ for some smooth conic $F$.
    When $C^2 = 11$, there exists a birational morphism $\pi: \underline{X} \to \PP^2$ such that $C$ is a member of $|-K_{\underline{X}} + \pi^*H'|$ where $H'$ is a hyperplane class of $\PP^2$.
\end{proposition}

\begin{proof}
    \begin{enumerate}
        \item \textit{Case 1: $C^2 = 5$.} In this case $2C$ is big and nef. Using Kawamata-Viehweg vanishing and Riemann-Roch we can conclude that $2C + K_{\underline{X}}$ is effective since $h^0(\underline{X}, \OO_{\underline{X}}(2C + K_{\underline{X}})) = 4 > 0$. We also have $$C \cdot (2C + K_{\underline{X}}) = 3, \quad -K_{\underline{X}} \cdot (2C + K_{\underline{X}}) = 10, \quad (2C + K_{\underline{X}})^2 = -4.$$
        After writing the Zariski decomposition of the divisor $2C + K_{\underline{X}}$, by similar calculations as in proposition \ref{dp3d5}, Case 1, we can conclude that the number of $N_i$'s is $4$.
        This is only possible when $\underline{X}'$ is blow-up of $\PP^2$ at one point. Let $H'$ be the pullback of a hyperplane class in $\PP^2$ via a blow up $\beta':\underline{X}' \to \PP^2$ at one point and $E$ be the unique exceptional divisor of $\beta'$, then $C' \sim 3H' - 2E$ and $C \in |3\pi^*H' - 2\pi^*E|$.
        For different contraction maps, they are one of the following forms. Using \ref{prog} we get:
        \begin{enumerate}
            \item $3H - 2E_i$, where $1 \leq i \leq 5$.
            \item $4H - 3E_i - E_j - E_k$, where $1 \leq i, j, k \leq 5$ are distinct indices.
            \item $5H - 4E_i - E_j - E_k - E_l - E_m$.
            \item $6H - 3E_i - 3E_j - 3E_k - 2E_l$, where $1 \leq i, j, k, l \leq 5$ are distinct indices.
            \item $7H - 4E_i - 3E_j - 3E_k - 3E_l - E_m$,  where $1 \leq i, j, k, l, m \leq 5$ are distinct indices.
        \end{enumerate}
        \item \textit{Case 2: $C^2 = 7$.} In this case $C+K_{\underline{X}}$ is disjoint union of $3$ $(-1)$-curves. Let $\pi: \underline{X} \to \underline{X}'$ be a birational contraction of these $(-1)$-curves. Then $\underline{X}'$ is a del Pezzo surface of degree $7$ and $C$ is a member of $|-\pi^* K_{\underline{X}'}|$. For different contraction maps, we get $C$ as one of the following forms. Using \ref{prog} we get:
        \begin{enumerate}
            \item $3H - E_i - E_j$,  where $1 \leq i, j \leq 5, i \neq j$.
            \item $4H - 2E_i - 2E_j - E_k$,  where $1 \leq i, j, k \leq 5$ are distinct indices.
            \item $5H - 3E_i - 2E_j - 2E_k - E_l$,  where $1 \leq i, j, k, l \leq 5$ are distinct indices.
            \item $6H - 3E_i - 3E_j - 3E_k - E_l - E_m$, $6H - 4E_i - 2E_j - 2E_k - 2E_l - E_m$.
            \item $7H - 4E_i - 3E_j - 3E_k - 2E_l - 2E_m$,  where $1 \leq i, j, k, l, m \leq 5$ are distinct indices.
        \end{enumerate}
        \item \textit{Case 3: $C^2 = 9$.} In this case $C$ is big and nef. Using Kawamata-Viehweg vanishing and Riemann-Roch $C + K_{\underline{X}}$ is effective since $h^0(\underline{X}, \OO_{\underline{X}}(C + K_{\underline{X}})) = 2 > 0$. We also have $$C \cdot (C + K_{\underline{X}}) = 2, \quad -K_{\underline{X}} \cdot (C + K_{\underline{X}}) = 3, \quad (C + K_{\underline{X}})^2 = -1.$$
        After writing the Zariski decomposition of the divisor $C + K_{\underline{X}}$, by similar calculations as in proposition \ref{dp3d5}, Case 1, we can conclude that the number of $N_i$'s is one, say $N$. Let $\pi: \underline{X} \to \underline{X}'$ be a birational contraction of $N$ by $C$.
        This is only possible when $\underline{X}'$ is blow-up of $\PP^2$ at $4$ points. Let $C'$ be the pushforward of $C$ in $\underline{X}'$. Then using similar discussions as in proposition \ref{dp2d4}, case 3, we can conclude that $C'$ is linearly equivalent to $-K_{\underline{X}'} + F$ for some smooth conic $F$. Hence $C$ is a member of $|-\pi^*K_{\underline{X}'} + \pi^*F|$.
        For different contraction maps, they are one of the following forms. Using \ref{prog} we get:
        \begin{enumerate}
            \item $4H - 2E_i - E_j - E_k - E_l$, where $1 \leq i, j, k, l \leq 5$ are distinct indices.
            \item $5H - 2E_i - 2E_j - 2E_k - 2E_l$, $5H - 3E_i - 2E_j - E_k - E_l - E_m$.
            \item $6H - 3E_i - 3E_j - 2E_k - 2E_l - E_m$.
            \item $7H - 3E_i - 3E_j - 3E_k - 3E_l - 2E_m$,  where $1 \leq i, j, k, l, m \leq 5$ are distinct indices.
        \end{enumerate}
        \item \textit{Case 4: $C^2 = 11.$} In this case we have $$(C+K_{\underline{X}})^2 = 1, \quad -K_{\underline{X}} \cdot (C + K_{\underline{X}}) = 3.$$ Therefore $C+K_{\underline{X}}$ is linearly equivalent to pullback of a hyperplane class $H'$ of $\PP^2$ via a birational morphism $\pi: \underline{X} \to \PP^2$. That is $C \sim -K_{\underline{X}} + \pi^* H'$. All such curves are one of the following forms. Using \ref{prog} we get:
        \begin{enumerate}
            \item $4H - E_i - E_j - E_k - E_l - E_m$.
            \item $5H - 2E_i - 2E_j - 2E_k - E_l - E_m$.
            \item $6H - 3E_i - 2E_j - 2E_k - 2E_l - 2E_m$, where $1 \leq i, j, k, l, m \leq 5$ are distinct indices.
        \end{enumerate}
    \end{enumerate}
\end{proof}

Next assume that $C$ is an irreducible curve of anticanonical degree $8$. Then by Hodge Index theorem and adjunction formula, we must have $C^2 = 6,8,10,12,14$ or $16$.

\begin{proposition}
    \label{dp4d8}
    Let $\underline{X}$ be a del Pezzo surface of degree $4$ and $C$ is an irreducible curve of anticanonical degree $8$. Then we have either $C^2 = 6,8,10,12,14$ or $16.$
    When $C^2 = 6$ there exists a birational morphism $\pi: \underline{X} \to \underline{X}'$ to $\PP^1 \times \PP^1$ such that $C$ is in pullback of a member in the class $(3,1)$ or $(1,3)$ in $\PP^1 \times \PP^1$.
    When $C^2 = 8$ there exists a birational morphism $\pi: \underline{X} \to \underline{X}'$ to blow up of $\PP^2$ at a point such that $C$ is a member of $|-\pi^*K_{\underline{X}'}|$.
    When $C^2 = 10$ there exists a birational morphism $\pi: \underline{X} \to \underline{X}'$ to blow up of $\PP^2$ at three points such that $C$ is a member of $|-\pi^*K_{\underline{X}'} + \pi^*F|$ for some smooth conic $F$.
    When $C^2 = 12$, then one can find a smooth conic $F$ such that $C \sim -K_{\underline{X}} + 2F$.
    When $C^2 = 14$ there exists a birational map $\pi: \underline{X} \to \PP^1 \times \PP^1$ such that $C$ is a member of $|-K_{\underline{X}} + \pi^*(1,1)|$.
    When $C^2 = 16$, $C \in |-2K_{\underline{X}}|$.
\end{proposition}

\begin{proof}
    \begin{enumerate}
        \item \textit{Case 1: $C^2 = 6$.} Using similar analysis as in proposition \ref{dp3d6}, case 1, we can conclude that $C$ is pullback of a member either in $(3,1)$-class or in $(1,3)$-class in $\PP^1 \times \PP^1$ via a birational morphism. 
        They are one of the following forms. Using \ref{prog} we get:
        \begin{enumerate}
            \item $4H - 3E_i - E_j$,  where $1 \leq i, j \leq 5, i \neq j$.
            \item $5H - 4E_i - E_j - E_k - E_l$.
            \item $7H - 4E_i - 3E_j - 3E_k - 3E_l$,  where $1 \leq i, j, k, l \leq 5$ are distinct indices.
            \item $8H - 4E_i - 4E_j - 4E_k - 3E_l - E_m$, where $1 \leq i, j, k, l, m \leq 5$ are distinct indices.
        \end{enumerate}
        \item \textit{Case 2: $C^2 = 8$.} In this case $C+K_{\underline{X}}$ is union of $4$ disjoint $(-1)$-curves. Let $\pi: \underline{X} \to \underline{X}'$ be a birational contraction of these curves. Then $\underline{X}'$ is a del Pezzo surface of degree $8$ and $C$ is a member of $|-\pi^*K_{\underline{X}'}|$.
        For different contractions we get the following type of curves. Using \ref{prog} we get:
        \begin{enumerate}
            \item $3H - E_i$,  where $1 \leq i \leq 5$.
            \item $4H - 2E_i - 2E_j$, where $1 \leq i, j \leq 5, i \neq j$.
            \item $5H - 3E_i - 2E_j - 2E_k$,  where $1 \leq i, j, k \leq 5$ are distinct indices.
            \item $6H - 3E_i - 3E_j - 3E_k - E_l;$ $6H - 4E_i - 2E_j - 2E_k - 2E_l;$ $6H - 4E_i - 3E_j - E_k - E_l - E_m$.
            \item $7H - 5E_i - 2E_j - 2E_k - 2E_l - 2E_m$.
            \item $8H - 4E_i - 4E_j - 4E_k - 2E_l - 2E_m;$ $8H - 5E_i - 3E_j - 3E_k - 3E_l - 2E_m$, where $1 \leq i, j, k, l, m \leq 5$ are distinct indices.
        \end{enumerate}
        \item \textit{Case 3: $C^2 = 10$.} In this case $C$ is big and nef. Using Kawamata-Viehweg vanishing and Riemann-Roch we conclude that $C+K_{\underline{X}}$ is effective since $h^0(\underline{X}, \OO_{\underline{X}}(C + K_{\underline{X}})) = 2 > 0.$ We also have $$C \cdot (C + K_{\underline{X}}) = 2, \quad -K_{\underline{X}} \cdot (C + K_{\underline{X}}) = 4, \quad (C + K_{\underline{X}})^2 = -2.$$
        After writing the Zariski decomposition of the divisor $C + K_{\underline{X}}$, by similar calculations as in proposition \ref{dp3d5}, Case 1, we can conclude that the number of $N_i$'s is $2$. Let $\pi: \underline{X} \to \underline{X}'$ be a birational contraction of $N_i$'s by $C$.
        This is only possible when $\underline{X}'$ is blow-up of $\PP^2$ at $3$ points. Let $C'$ be the pushforward of $C$ in $\underline{X}'$. Using similar discussions as in proposition \ref{dp2d4}, case 3, we can conclude that $C'$ is linearly equivalent to $-K_{\underline{X}'} + F$ for some smooth conic $F$. Therefore $C$ is a member of $|-\pi^*K_{\underline{X}'} + \pi^*F|$ for some smooth conic $F$. 
        They are one of the following forms. Using \ref{prog} we get:
        \begin{enumerate}
            \item $4H - 2E_i - E_j - E_k$,  where $1 \leq i, j, k \leq 5$ are distinct indices.
            \item $5H - 3E_i - 2E_j - E_k - E_l$,  where $1 \leq i, j, k, l \leq 5$ are distinct indices.
            \item $6H - 3E_i - 3E_j - 2E_k - 2E_l;$ $6H - 4E_i - 2E_j - 2E_k - E_l - E_m$.
            \item $7H - 4E_i - 3E_j - 3E_k - 2E_l - E_m$.
            \item $8H - 4E_i - 4E_j - 3E_k - 3E_l - 2E_m$, where $1 \leq i, j, k, l, m \leq 5$ are distinct indices.
        \end{enumerate}
        \item \textit{Case 4: $C^2 = 12$.} Let $D = \frac{C + K_{\underline{X}}}{2}$ be a $\Q$-divisor. Then $2D$ is an integral divisor and since in this case $C$ is big and nef, using Kawamata-Viehweg vanishing and Riemann-Roch we can conclude that $2D$ is effective. Moreover, $h^0(\underline{X}, \OO_{\underline{X}}(2D)) = 3 > 0$. Hence $\lvert 2D \rvert$ defines a birational morphism $\pi: \underline{X} \to \PP^2$ and $\kappa(\OO_{\underline{X}}(2D)) \geq 1$, where $\kappa$ is the Iitaka dimension of the line bundle $\OO_{\underline{X}}(2D)$ corresponding to the divisor $2D$. We also have $D^2 = 0$ and $-K_{\underline{X}} \cdot D = 2$. Hence $\kappa(\OO_{\underline{X}}(2D)) \neq 2$, therefore $\kappa(\OO_{\underline{X}}(2D)) = 1$. Hence we get a proper morphism $\pi: \underline{X} \to \PP^2$ such that $\pi_{\lvert 2D \rvert}(\underline{X})$ has dimension $1$. So we have a Stein factorization $$\underline{X} \xrightarrow{g} B \xrightarrow{h} \pi_{\lvert 2D \rvert} (\underline{X}) \subset \PP^2$$ where $g$ has geometrically connected fibers, $B$ is a smooth projective curve and $h$ is finite.
        Since $\underline{X}$ is a del Pezzo surface and $g$ has connected fibers, $B$ must be $\PP^1$. Hence we have $$\underline{X} \xrightarrow{g} \PP^1 \xrightarrow{h} \pi_{\lvert 2D \rvert} (\underline{X}).$$ Now, since $g$ has geometrically connected fibers and $\kappa(\OO_{\underline{X}}(2D)) = 1$, we have $g_*\OO_{\underline{X}} = \OO_{\PP^1}$ and $\OO_{\underline{X}}(2D) = g^*L$ for some line bundle $L$ on $\PP^1$. 
        Hence $$h^0(\underline{X}, \OO_{\underline{X}}(2D)) = h^0(\PP^1, g_*\OO_{\underline{X}}(2D)) = h^0(\PP^1, L) = 3.$$ Therefore $L = \OO_{\PP^1}(2)$. Hence $\OO_{\underline{X}}(2D) = g^*\OO_{\PP^1}(2)$ and $D$ is linearly equivalent to $F$, where $F$ is a smooth fiber of the conic fibration $g$. So, $C \sim -K_{\underline{X}} + 2F$ for some smooth conic $F$. 
        Using \ref{prog} we get the following forms:
        \begin{enumerate}
            \item $4H - E_i - E_j - E_k - E_l$.
            \item $5H - 2E_i - 2E_j - 2E_k - E_l$, where $1 \leq i, j, k, l \leq 5$ are distinct indices.
            \item $6H - 3E_i - 3E_j - 2E_k - E_l - E_m$.
            \item $7H - 4E_i - 3E_j - 2E_k - 2E_l - 2E_m$.
            \item $5H - 3E_i - E_j - E_k - E_l - E_m$.
            \item $7H - 3E_i - 3E_j - 3E_k - 3E_l - E_m$, where $1 \leq i, j, k, l, m \leq 5$ are distinct indices.
        \end{enumerate}
        \item \textit{Case 5: $C^2 = 14$.} In this case $(C+K_{\underline{X}})^2 = 2$ and $-K_{\underline{X}} \cdot (C+K_{\underline{X}}) = 4$. Hence there exists a birational map $\pi: \underline{X} \to \PP^1 \times \PP^1$ such that $C + K_{\underline{X}}$ is pullback of a member in $(1,1)$ class in $\PP^1 \times \PP^1$. Hence $C \in |-K_{\underline{X}} + \pi^*(1,1)|$. For different pullbacks we get the following forms. Using \ref{prog} we get:
        \begin{enumerate}
            \item $5H - 2E_i - 2E_j - E_k - E_l - E_m$.
            \item $6H - 3E_i - 2E_j - 2E_k - 2E_l - E_m$.
            \item $7H - 3E_i - 3E_j - 3E_k - 2E_l - 2E_m$, where $1 \leq i, j, k, l, m \leq 5$ are distinct indices.
        \end{enumerate}
        \item \textit{Case 6: $C^2 = 16$.} In this case $C \in |-2K_{\underline{X}}|$.
    \end{enumerate}
\end{proof}

\subsection{Classification of irreducible boundaries.} Next we let $\Delta = \epsilon D$ where $D$ is an irreducible curve on $\underline{X}$ and $\epsilon = 1 - 1/m$. We discuss whether $-(K_{\underline{X}} + \Delta)$ satisfies some positivity of divisors when $D$ is an irreducible curve of anticanonical degree at most $8$.

Using Lemma \ref{thelemma} for a del Pezzo Surface of degree $4$, we can conclude that a divisor $L$ is ample (nef) if and only if $L \cdot E > 0$ (resp. $L \cdot E \geq 0$) for any $(-1)$-curve $E$. 
We use \ref{prog} to calculate the intersection numbers.

\begin{itemize}
    \item \textit{Case of $-K_{\underline{X}} \cdot D = 1:$} We assume that $D$ is a $(-1)$-curve. In this case $D^2 = -1$. Since $D \cdot E = 0, 1$ for any $(-1)$-curve $E$ distinct from $D$, we conclude that $-(K_{\underline{X}} + \epsilon D)$ is ample for any $\epsilon = 1 - 1/m$ where $m \in \Z_{\geq 1}.$ In this case $-(K_{\underline{X}} + D)$ is big and nef but not ample.
    \item \textit{Case of $-K_{\underline{X}} \cdot D = 2:$} We assume that $D$ is a smooth fiber of a conic fibration. In this case $D \cdot E = 0, 1$ for any $(-1)$-curve $E$. Hence we conclude that $-(K_{\underline{X}} + \epsilon D)$ is ample for any $\epsilon = 1 - 1/m$ where $m \in \Z_{\geq 1}$.
    \item \textit{Case of $-K_{\underline{X}} \cdot D = 3:$} We assume that $D$ is pullback of a hyperplane class in $\PP^2$ via a birational morphism. In this case $D \cdot E = 0, 1$ and $2$ for any $(-1)$-curve $E$. Hence we conclude that $-(K_{\underline{X}} + \epsilon D)$ is big and nef for $\epsilon = \frac{1}{2}.$
    \item \textit{Case of $-K_{\underline{X}} \cdot D = 4:$} We first assume that $D$ is pullback of a member in $(1,1)$ class of $\PP^1 \times \PP^1$. In this case $D \cdot E = 0, 1$ and $2$ for any $(-1)$-curve $E$. Hence we conclude that $-(K_{\underline{X}} + \epsilon D)$ is big and nef for $\epsilon = \frac{1}{2}.$

    Now we assume that $D \in |-K_{\underline{X}}|$, since $D \cdot E = 1$ for any $(-1)$-curve $E$, we conclude that $-(K_{\underline{X}} + \epsilon D)$ is ample for any $\epsilon = 1 - 1/m$ where $m \in \Z_{\geq 1}.$
    \item \textit{Case of $-K_{\underline{X}} \cdot D = 5:$} We assume that $D^2 = 3$, then $D$ is the pullback of a divisor of the form $2H' - E$ as explained above, where $H'$ is the pullback of a hyperplane class from $\PP^2$ via a blow up $\beta':\underline{X}' \to \PP^2$ at one point and $E$ is the exceptional divisor of $\beta'$. For a $(-1)$-curve $F$, $D \cdot F = 0, 1, 2$ or $3$. So we conclude that $-(K_{\underline{X}} + \epsilon D)$ is not nef for any $\epsilon > 0.$

    Now we assume that $D$ is pullback of a member of the anticanonical class in a del Pezzo surface of degree $5$. Then $D \cdot E = 0, 1$ or $2$. Hence, $-(K_{\underline{X}} + \epsilon D)$ is nef and big for $\epsilon = \frac{1}{2}.$
    \item \textit{Case of $-K_{\underline{X}} \cdot D = 6:$} We first assume that $D^2 = 4$. In this case, $D \cdot E = 0, 1, 2, 3$ and $4$ for any $(-1)$-curve $E$.
    In this case, $-(K_{\underline{X}} + \epsilon D)$ is not nef for any $\epsilon > 0$.

    Next we assume that $D$ is pullback of a member in the anticanonical class of a smooth sextic surface. In this case $D \cdot E = 0, 1, 2$ and $3$ for any $(-1)$-curve $E$. So we conclude that $-(K_{\underline{X}} + \epsilon D)$ is not nef for any $\epsilon > 0.$

    Next we assume that $D$ is linearly equivalent to $-K_{\underline{X}} + F$ for some smooth conic $F$. In this case $D$ is ample and $D \cdot E = 1$ and $2$ for any $(-1)$-curve $E$. Hence we can conclude that $-(K_{\underline{X}} + \epsilon D)$ is nef for $\epsilon = \frac{1}{2}$, but it is not big.
    \item \textit{Case of $-K_{\underline{X}} \cdot D = 7:$} We first assume that $D^2 = 5$. Then $D$ is pullback of a divisor of the form $3H' - 2E$ as described above, where $H'$ is the pullback of a hyperplane class of $\PP^2$ via a blow up $\beta':\underline{X}' \to \PP^2$ at one point and $E$ is the exceptional divisor of $\beta'$. In this case, $D \cdot F = 0 ,1, 2, 3$ and $4$ for any $(-1)$-curve $F$. Hence $-(K_{\underline{X}} + \epsilon D)$ is not nef for any $\epsilon > 0.$

    Next we assume that $D$ is a member of the anticanonical class of a del Pezzo surface of degree $7$. In this case $D \cdot E = 0, 1, 2, 3$ and $4$ for any $(-1)$-curve $E$. Therefore we can conclude that $-(K_{\underline{X}} + \epsilon D)$ is not nef for any $\epsilon > 0.$

    Next we assume that $D^2 = 9$. In this case $D \cdot E = 0, 1, 2$ and $3$ for any $(-1)$-curve $E$. Hence we conclude that $-(K_{\underline{X}} + \epsilon D)$ is not nef for any $\epsilon > 0.$

    Next we assume that $D^2 = 11$. In this case $D + K_{\underline{X}}$ is linearly equivalent to pullback of a hyperplane class $H$ of $\PP^2$ as described above. In this case $D \cdot E = 1, 2$ and $3$. We conclude that $-(K_{\underline{X}} + \epsilon D)$ is not nef for any $\epsilon > 0$.
    \item \textit{Case of $-K_{\underline{X}} \cdot D = 8:$} We assume that $D^2 = 6$. In this case $D$ is pullback of a divisor in the class $(3,1)$ or $(1,3)$ via a birational morphism to $\PP^1 \times \PP^1$ as described above. In this case $D \cdot F = 0, 1, 3$ and $4$ for any $(-1)$-curve $F$. So we conclude that $-(K_{\underline{X}} + \epsilon D)$ is not nef for any $\epsilon > 0.$

    Next we assume that $D$ is a member of the anticanonical class of a del Pezzo surface of degree $8$. In this case $D \cdot E$ can be as high as $5$. So $-(K_{\underline{X}} + \epsilon D)$ is not nef for any $\epsilon > 0.$

    Next we assume that $D^2 = 10$. In this case $D \cdot E = 0, 1, 2, 3$ and $4$ for any $(-1)$-curve $E$. Hence we conclude that $-(K_{\underline{X}} + \epsilon D)$ is not nef for any $\epsilon > 0.$

    Next we assume that $D^2 = 12$. Then $D$ is linearly equivalent to $-K_{\underline{X}} + 2F$ for some smooth conic $F$. In this case $D \cdot E = 1, 3$ for two choices of the conic and $D \cdot E = 0, 1, 2, 3$ and $4$ for other $4$ choices, for any $(-1)$-curve $E$. Hence we conclude that in any case $-(K_{\underline{X}} + \epsilon D)$ is not nef for any $\epsilon > 0.$

    Next we assume that $D^2 = 14$. In this case $D+K_{\underline{X}}$ is pullback of the $(1,1)$ class in $\PP^1 \times \PP^1$ as described above. In this case $D \cdot E = 1, 2$ and $3$ for any $(-1)$-curve $E$. Hence we conclude that $-(K_{\underline{X}} + \epsilon D)$ is not nef for any $\epsilon > 0$.

    Next we assume that $D \in |-2K_{\underline{X}}|$. Since $D \cdot E = 2$ for any $(-1)$-curve $E$, we conclude that $-(K_{\underline{X}} + \epsilon D)$ is nef for $\epsilon = \frac{1}{2}$, but it is not big.
 \end{itemize}

In summary, we conclude:

\begin{proposition}
    Let $X$ is a del Pezzo orbifold whose underlying surface $\underline{X}$ is a del Pezzo surface of degree $4$. Let $D$ be an irreducible curve of degree less than or equal to $8$, and $\epsilon = 1 - 1/m, m \in \Z_{\geq 2}$. Then $-(K_{\underline{X}} + \epsilon D)$ is ample if and only if
    \begin{itemize}
        \item $D$ is a $(-1)$-curve and $\epsilon > 0$.
        \item $D$ is a smooth fiber of a conic fibration and $\epsilon > 0$.
        \item $D$ is a member of the anticanonical class and $\epsilon > 0$.
    \end{itemize}

    Also, $-(K_{\underline{X}} + \epsilon D)$ is nef only when $\epsilon = \frac{1}{2}$ and $D$ is one of the following.

    \begin{itemize}
        \item pullback of a hyperplane class in $\PP^2$. $-(K_{\underline{X}} + \epsilon D)$ is big in this case.
        \item pullback of a member in $(1,1)$-class of $\PP^1 \times \PP^1$. $-(K_{\underline{X}} + \epsilon D)$ is big in this case.
        \item linearly equivalent to $-K_{\underline{X}} + F$ for some smooth conic $F$, $-(K_{\underline{X}} + \epsilon D)$ is nef but not big.
        \item member of $|-2K_{\underline{X}}|$, $-(K_{\underline{X}} + \epsilon D)$ is nef but not big.
    \end{itemize}
\end{proposition}

\section{Surfaces of degree 5.}
\label{dp5}

Suppose $\underline{X}$ is a del Pezzo surface of degree $5$ and $\Delta_{\epsilon}$ be a $\Q$-divisor on it. In this case $\underline{X}$ is the blow up $\beta: \underline{X} \to \PP^2$ at four general points. The canonical divisor is given by $$K_{\underline{X}} = -3H + E_1 + E_2 + E_3 + E_4,$$ where
$E_i$'s $(1 \leq i \leq 4)$ are exceptional curves of $\beta$ and $\beta^*K_{\PP^2} = -3H$, where $H$ is the hyperplane class in $\underline{X}$. We classify all such del Pezzo orbifolds $(\underline{X}, \Delta_{\epsilon})$. $\underline{X}$ can be realised as blow up of $\PP^2$ at $4$ general points.
In this case $K_{\underline{X}}$ is of degree $5$, so $\Delta_{\epsilon}$ can have at most $10$ components and each component has anticanonical degree at most $10$.

\subsection{Low degree curves on $\underline{X}$.} We first classify irreducible curves on $\underline{X}$ with anticanonical at most $4$.

\begin{itemize}
    \item an irreducible curve $C$ of anticanonical degree $1$ is a $(-1)$-curve such that $C^2 = -1$. There are $10$ $(-1)$-curves which are of the following form. Using \ref{prog} we get:
    \begin{enumerate}
        \item $E_i$, which are $4$ exceptional curves under the blow up $\beta: \underline{X} \to \PP^2$ at four general points.
        \item $H - E_i - E_j$, where $1 \leq i, j \leq 4, i \neq j$. There are $6$ such curves.
    \end{enumerate}
    \item an irreducible curve $C$ of anticanonical degree $2$ must satisfy $C^2 = 0$ by Hodge Index theorem and adjunction formula. In this case $|C|$ is a basepoint-free pencil. Hence it defines a morphism $f: \underline{X} \to \PP^1$. Such $C$ is a smooth fiber of a conic fibration $f$ on $\underline{X}$.
    They are one of the following forms. Using \ref{prog} we get:
    \begin{enumerate}
        \item $H - E_i$, where $1 \leq i \leq 4$.
        \item $2H - E_1 - E_2 - E_3 - E_4$.
    \end{enumerate}
    \item an irreducible curve $C$ of anticanonical degree $3$ must satisfy $C^2 = 1$ by Hodge Index theorem and adjunction formula. 
    Similarly as in proposition \ref{dp2d3}, Case 1, we can conclude that $C$ is a member of the pullback of a hyperplane class in $\PP^2$ via a birational morphism.
    For different birational contractions they are one of the following forms. Using \ref{prog} we get:
    \begin{enumerate}
        \item $H$.
        \item $2H - E_i - E_j - E_k$, where $1 \leq i, j, k \leq 4$ are distinct indices.
    \end{enumerate}
    \item an irreducible curve $C$ of anticanonical degree $4$ must satisfy $C^2 = 2$ by Hodge Index theorem and adjunction formula. 
    In this case, similarly as in proposition \ref{dp2d4}, Case 1, we can conclude that $C$ is pullback of a member of the $(1,1)$-class in $\PP^1 \times \PP^1$ via a birational morphism.
    For different pullbacks they can be described as one of the following forms. Using \ref{prog} we get:
    \begin{enumerate}
        \item $2H - E_i - E_j$, where $1 \leq i, j \leq 4, i \neq j$.
        \item $3H - 2E_i - E_j - E_k - E_l$, where $1 \leq i, j, k, l \leq 4$ are distinct indices.
    \end{enumerate}    
\end{itemize}

Next assume that $C$ is an irreducible curve of anticanonical degree $5$. Then it must satisfy $C^2 = 3$ or $C^2 = 5$ by Hodge Index theorem and adjunction formula.

\begin{proposition}
    Let $\underline{X}$ be a del Pezzo surface of degree $5$ and $C$ be an irreducible curve of anticanonical degree $5$. Then $C^2 = 3$ or $5$. In the former case, there exists a birational morphism $\pi: \underline{X} \to \underline{X}'$ to blow up of $\PP^2$ at one point and $C$ is a member of $|2\pi^*H' - \pi^* E|$ where $H'$ is the pullback of the hyperplane class on $\PP^2$ via a blow up $\beta':\underline{X}' \to \PP^2$ at one point and $E$ is the unique exceptional divisor of $\beta'$,
    and in the later case $C$ is a member of the anticanonical class.
\end{proposition}

\begin{proof}
    \begin{enumerate}
        \item \textit{Case 1: $C^2 = 3$.} This is similar to the proof of proposition \ref{dp3d5}. Using \ref{prog} we get the following forms:
        \begin{enumerate}
            \item $2H - E_i$, where $1 \leq i \leq 4$.
            \item $3H - 2E_i - E_j - E_k$, where $1 \leq i, j, k \leq 4$ are distinct indices.
            \item $4H - 2E_i - 2E_j - 2E_k - E_l$, where $1 \leq i, j, k, l \leq 4$ are distinct indices.
        \end{enumerate}
        \item \textit{Case 2: $C^2 = 5$.} In this case $C \in |-K_{\underline{X}}|$.
    \end{enumerate}
\end{proof}

Next we assume that $C$ is an irreducible curve of anticanonical degree $6$. Then $C^2 = 4$ or $C^2 = 6$ by Hodge Index theorem and adjunction formula.

\begin{proposition}
    Let $\underline{X}$ be a del Pezzo surface of degree $5$ and $C$ is an irreducible curve of anticanonical degree $6$. Then either we have $C^2 = 4$ or $C^2 = 6.$ 
    In the former case, there exists a birational morphism $\pi: \underline{X} \to \underline{X}'$ to $\PP^1 \times \PP^1$ such that $C$ is a pullback of a member of $(2,1)$ or $(1,2)$-class in $\PP^1 \times \PP^1$. 
    and in the later case, there exists a birational morphism $\pi: \underline{X} \to \underline{X}'$ to a smooth sextic del Pezzo surface such that $C$ is a member of $|-\pi^*K_{\underline{X}'}|$.
\end{proposition}

\begin{proof}
    \begin{enumerate}
        \item \textit{Case 1: $C^2 = 4$.} This case is similar to the proof of proposition \ref{dp3d6}, case 1. Using \ref{prog} we get the following forms:
        \begin{enumerate}
            \item $2H$.
            \item $3H - 2E_i - E_j$, where $1 \leq i, j \leq 4, i \neq j$.
            \item $4H - 2E_i - 2E_j - 2E_k;$ $4H - 3E_i - E_j - E_k - E_l$.
            \item $5H - 3E_i - 2E_j - 2E_k - 2E_l$, where $1 \leq i, j, k, l \leq 4$ are distinct indices.
        \end{enumerate}
        \item \textit{Case 2: $C^2 = 6$.} This case is similar to the proof of proposition \ref{dp4d6}, case 2. Using \ref{prog} we get the following forms:
        \begin{enumerate}
            \item $3H - E_i - E_j - E_k$, where $1 \leq i, j, k \leq 4$ are distinct indices.
            \item $4H - 2E_i - 2E_j - E_k - E_l$, where $1 \leq i, j, k, l \leq 4$ are distinct indices.
        \end{enumerate}
    \end{enumerate}
\end{proof}

Next we assume that $C$ is an irreducible curve of anticanonical degree $7$. Then $C^2 = 5, 7$ or $9$ by Hodge Index theorem and adjunction formula.

\begin{proposition}
    Let $\underline{X}$ be a del Pezzo surface of degree $5$ and $C$ is an irreducible curve of anticanonical degree $7$. Then either we have $C^2 = 5, 7$ or $9$.
    When $C^2 = 5$, there exists a birational morphism $\pi: \underline{X} \to \underline{X}'$ to blow up of $\PP^2$ at one point such that $C$ is a member of $|3\pi^*H' - 2\pi^*E|$, where $H'$ is the pullback of a hyperplane class in $\PP^2$ via a blow up $\beta':\underline{X}' \to \PP^2$ at one point and $E$ is the unique exceptional divisor of $\beta'$.
    When $C^2 = 7$, there exists a birational morphism $\pi: \underline{X} \to \underline{X}'$ to a del Pezzo surface of degree $7$ such that $C$ is a member of $|-\pi^*K_{\underline{X}'}|$.
    When $C^2 = 9$, then $C$ is ample and one can find a smooth conic $F$ such that $C \sim -K_{\underline{X}} + F.$
\end{proposition}

\begin{proof}
    \begin{enumerate}
        \item \textit{Case 1: $C^2 = 5$.} This case is similar to the proof of the proposition \ref{dp4d7}. Using \ref{prog} we get the following forms:
        \begin{enumerate}
            \item $3H - 2E_i$, where $1 \leq i \leq 4$.
            \item $4H - 3E_i - E_j - E_k$, where $1 \leq i, j, k \leq 4$ are distinct indices.
            \item $6H - 3E_i - 3E_j - 3E_k - 2E_l$, where $1 \leq i, j, k, l \leq 4$ are distinct indices.
        \end{enumerate}
        \item \textit{Case 2: $C^2 = 7$.}  This case is similar to the proof of the proposition \ref{dp4d7}. Using \ref{prog} we get the following forms:
        \begin{enumerate}
            \item $3H - E_i - E_j$, where $1 \leq i, j \leq 4, i \neq j$.
            \item $4H - 2E_i - 2E_j - E_k$, where $1 \leq i, j, k \leq 4$ are distinct indices.
            \item $5H - 3E_i - 2E_j - 2E_k - E_l$, where $1 \leq i, j, k, l \leq 4$ are distinct indices.
        \end{enumerate}
        \item \textit{Case 3: $C^2 = 9$.} This case is similar to the proof of the proposition \ref{dp4d7}. Using \ref{prog} we get the following forms:
        \begin{enumerate}
            \item $4H - 2E_i - E_j - E_k - E_l$.
            \item $5H - 2E_i - 2E_j - 2E_k - 2E_l$, where $1 \leq i, j, k, l \leq 4$ are distinct indices.
        \end{enumerate}
    \end{enumerate}
\end{proof}

Next we assume that $C$ is an irreducible curve of anticanonical degree $8$. Then $C^2 = 6, 8, 10$ or $12$ by Hodge Index theorem and adjunction formula.

\begin{proposition}
    Let $\underline{X}$ be a del Pezzo surface of degree $5$ and $C$ is an irreducible curve of anticanonical degree $8$. Then $C^2 = 6, 8, 10$ or $12$.
    When $C^2 = 6$, there exists a birational morphism $\pi: \underline{X} \to \underline{X}'$ to $\PP^1 \times \PP^1$ such that $C$ is pullback of a member of $(3,1)$ or $(1,3)$-class.
    When $C^2 = 8$, there exists a birational morphism $\pi: \underline{X} \to \underline{X}'$ to blow up of $\PP^2$ at one point such that $C$ is a member of $|-\pi^*K_{\underline{X}'}|$.
    When $C^2 = 10$, there exists a birational morphism $\pi: \underline{X} \to \underline{X}'$ to a smooth sextic del Pezzo surface such that $C$ is a member of $|-\pi^*K_{\underline{X}'} + \pi^*F|$, for some smooth conic $F$.
    When $C^2 = 12$, there exists a birational morphism $\pi: \underline{X} \to \PP^2$ such that $C \in |-K_{\underline{X}} + \pi^* H'|$ where $H'$ is the hyperplane class in $\PP^2.$
\end{proposition}

\begin{proof}
    \begin{enumerate}
        \item \textit{Case 1: $C^2 = 6$.} This case is similar to the proof of proposition \ref{dp4d8}. Using \ref{prog} we get the following forms:
        \begin{enumerate}
            \item $4H - 3E_i - E_j$, where $1 \leq i, j \leq 4, i \neq j$.
            \item $5H - 4E_i - E_j - E_k - E_l$.
            \item $7H - 4E_i - 3E_j - 3E_k - 3E_l$, where $1 \leq i, j, k, l \leq 4$ are distinct indices.
        \end{enumerate}
        \item \textit{Case 2: $C^2 = 8$.} This case is similar to the proof of proposition \ref{dp4d8}. Using \ref{prog} we get the following forms:
        \begin{enumerate}
            \item $3H - E_i$, where $1 \leq i \leq 4$.
            \item $4H - 2E_i - 2E_j$, where $1 \leq i, j \leq 4, i \neq j$.
            \item $5H - 3E_i - 2E_j - 2E_k$, where $1 \leq i, j, k \leq 4$ are distinct indices.
            \item $6H - 3E_i - 3E_j - 3E_k - E_l;$ $6H - 4E_i - 2E_j - 2E_k - 2E_l$, where $1 \leq i, j, k, l \leq 4$ are distinct indices.
        \end{enumerate}
        \item \textit{Case 3: $C^2 = 10$.} This case is similar to the proof of proposition \ref{dp4d8}. Using \ref{prog} we get the following forms:
        \begin{enumerate}
            \item $4H - 2E_i - E_j - E_k$, where $1 \leq i, j, k \leq 4$ are distinct indices.
            \item $5H - 3E_i - 2E_j - E_k - E_l$.
            \item $6H - 3E_i - 3E_j - 2E_k - 2E_l$, where $1 \leq i, j, k, l \leq 4$ are distinct indices.
        \end{enumerate}
        \item \textit{Case 4: $C^2 = 12$.} In this case $(C + K_{\underline{X}})^2 = 1$ and $-K_{\underline{X}} \cdot (C + K_{\underline{X}}) = 3$. Hence $C + K_{\underline{X}}$ is linearly equivalent to pullback of a line in $\PP^2$ via a birational morphism $\pi: \underline{X} \to \PP^2.$ So, $C \in |-K_{\underline{X}} + \pi^* H'|$ where $H'$ is the hyperplane class in $\PP^2.$ Then $C$ can be one of the following forms. Using \ref{prog} we get:
        \begin{enumerate}
            \item $4H - E_i - E_j - E_k - E_l$.
            \item $5H - 2E_i - 2E_j - 2E_k - E_l$, where $1 \leq i, j, k, l \leq 4$ are distinct indices.
        \end{enumerate}
    \end{enumerate}
\end{proof}

Next we assume that $C$ is an irreducible curve of anticanonical degree $9$. Then $C^2 = 7, 9, 11, 13$ or $15$ by Hodge Index theorem and adjunction formula.

\begin{proposition}
    Let $\underline{X}$ be a del Pezzo surface of degree $5$ and $C$ is an irreducible curve of anticanonical degree $9$. Then $C^2 = 7, 9, 11, 13$ or $15$.
    When $C^2 = 7$, there exists a birational morphism $\pi: \underline{X} \to \underline{X}'$ to blow up of $\PP^2$ at one point such that $C$ is a member of $|4\pi^*H' - 3\pi^*E|$ where $H'$ is the pullback of a hyperplane class in $\PP^2$ via a blow up $\beta':\underline{X}' \to \PP^2$ at one point and $E$ is the unique exceptional divisor of $\beta'$.
    When $C^2 = 9$, there exists a birational morphism $\pi: \underline{X} \to \underline{X}'$ to $\PP^2$ such that $C$ is a member of $|-\pi^*K_{\underline{X}'}|$.
    When $C^2 = 11$, there exists a birational morphism $\pi: \underline{X} \to \underline{X}'$ to blow up of $\PP^2$ at two points such that $C$ is a member of $|-\pi^*K_{\underline{X}'} + \pi^*F|$, for some smooth conic $F$.
    When $C^2 = 13$, then $C$ is ample and one can find a smooth conic $F$ such that $C \sim -K_{\underline{X}} + 2F$.
    When $C^2 = 15$, there exists a birational map $\pi: \underline{X} \to \PP^1 \times \PP^1$ such that $C$ is a member of $|-K_{\underline{X}} + \pi^*(1,1)|$.
\end{proposition}

\begin{proof}
    \begin{enumerate}
        \item \textit{Case 1: $C^2 = 7$.} In this case $2C$ is big and nef. Hence using Kawamata-Viehweg vanishing and Riemann-Roch we can conclude that $2C + K_{\underline{X}}$ is effective. We also have $$C \cdot (2C + K_{\underline{X}}) = 5, \quad -K_{\underline{X}} \cdot (2C + K_{\underline{X}}) = 13, \quad (2C + K_{\underline{X}})^2 = -3.$$
        After writing the Zariski decomposition of the divisor $2C + K_{\underline{X}}$, by similar calculations as in proposition \ref{dp3d5} we can conclude that the number of $N_i$'s is $3$. Let $\pi: \underline{X} \to \underline{X}'$ be a birational contraction of $N_i$'s.
        This is only possible when $\underline{X}'$ is blow-up of $\PP^2$ at one point. Let $H'$ be the pullback of a hyperplane class in $\PP^2$ via a blow up $\beta':\underline{X}' \to \PP^2$ at one point and $E$ be the unique exceptional divisor of $\beta'$, then $C' \sim 4H' - 3E$ and $C \in |4\pi^*H' - 3\pi^*E|$.
        For different contraction maps, they are one of the following forms. Using \ref{prog} we get:
        \begin{enumerate}
            \item $4H - 3E_i$, where $1 \leq i \leq 4$.
            \item $5H - 4E_i - E_j - E_k$, where $1 \leq i, j, k \leq 4$ are distinct indices.
            \item $8H - 4E_i - 4E_j - 4E_k - 3E_l$, where $1 \leq i, j, k, l \leq 4$ are distinct indices.
        \end{enumerate}
        \item \textit{Case 2: $C^2 = 9$.} In this case $C+K_{\underline{X}}$ is disjoint union of $4$ $(-1)$-curves. Let $\pi: \underline{X} \to \underline{X}'$ be a birational contraction of these $(-1)$-curves. Then $\underline{X}'$ must be $\PP^2$ and $C$ is a member of $|-\pi^*K_{\PP^2}|.$ For different contractions we get the following forms. Using \ref{prog} we get:
        \begin{enumerate}
            \item $3H$.
            \item $6H - 3E_i - 3E_j - 3E_k$, where $1 \leq i, j, k \leq 4$ are distinct indices.
        \end{enumerate}
        \item \textit{Case 3: $C^2 = 11$.} In this case $C$ is big and nef. Hence using Kawamata-Viehweg vanishing and Riemann-Roch we can conclude that $|C + K_{\underline{X}}|$ is a pencil. We also have $$C \cdot (C + K_{\underline{X}}) = 2, \quad -K_{\underline{X}} \cdot (C + K_{\underline{X}}) = 4, \quad (C + K_{\underline{X}})^2 = -2.$$
        After writing the Zariski decomposition of the divisor $C + K_{\underline{X}}$, by similar calculations as in proposition \ref{dp3d5} we can conclude that the number of $N_i$'s is $2$. Let $\pi: \underline{X} \to \underline{X}'$ be a birational contraction of $N_i$'s by $C$.
        This is only possible when $\underline{X}'$ is blow-up of $\PP^2$ at two points. Let $C'$ be the pushforward of $C$ in $\underline{X}'$. Then using similar discussions as in proposition \ref{dp2d4}, case 3, we can conclude that $C'$ is linearly equivalent to $-K_{\underline{X}'} + F$ for some smooth conic $F$. Therefore $C$ is a member of $|-\pi^*K_{\underline{X}'} + \pi^*F|$ for some smooth conic $F$.
        For different contraction maps, they are one of the following forms. Using \ref{prog} we get:
        \begin{enumerate}
            \item $4H - 2E_i - E_j$, where $1 \leq i, j \leq 4, i \neq j$.
            \item $5H - 3E_i - 2E_j - E_k$, where $1 \leq i, j, k \leq 4$ are distinct indices.
            \item $6H - 4E_i - 2E_j - 2E_k - E_l$.
            \item $7H - 4E_i - 3E_j - 3E_k - 2E_l$, where $1 \leq i, j, k, l \leq 4$ are distinct indices.
        \end{enumerate}
        \item \textit{Case 4: $C^2 = 13$.} This proof is similar to proposition \ref{dp4d8}, case 4. $C$ is one of the following forms. Using \ref{prog} we get:
        \begin{enumerate}
            \item $4H - E_i - E_j - E_k$, where $1 \leq i, j, k \leq 4$ are distinct indices.
            \item $5H - 2E_i - 2E_j - 2E_k;$ $5H - 3E_i - E_j - E_k - E_l$.
            \item $6H - 3E_i - 3E_j - 2E_k - E_l$.
            \item $7H - 3E_i - 3E_j - 3E_k - 3E_l$, where $1 \leq i, j, k, l \leq 4$ are distinct indices.
        \end{enumerate}
        \item \textit{Case 5: $C^2 = 15$.} In this case $C$ is big and nef. Using Kawamata-Viehweg vanishing and Riemann-Roch we can conclude that $C + K_{\underline{X}}$ is effective. We also have $$-K_{\underline{X}} \cdot (C + K_{\underline{X}}) = 4, \quad (C + K_{\underline{X}})^2 = 2.$$ Hence we can conclude there exists a birational map $\pi: \underline{X} \to \PP^1 \times \PP^1$ such that $C + K_{\underline{X}}$ is a member of the pullback of $(1,1)$-class. That is, $C$ is a member of $|-K_{\underline{X}} + \pi^*(1,1)|$. For different birational maps, we get the following forms of $C$. Using \ref{prog} we get:
        \begin{enumerate}
            \item $5H - 2E_i - 2E_j - E_k - E_l$.
            \item $6H - 3E_i - 2E_j - 2E_k - 2E_l$, where $1 \leq i, j, k, l \leq 4$ are distinct indices.
        \end{enumerate}
    \end{enumerate}
\end{proof}

Next we assume that $C$ is an irreducible curve of anticanonical degree $10$. Then $C^2 = 8, 10, 12, 14, 16, 18$ or $20$ by Hodge Index theorem and adjunction formula.

\begin{proposition}
    Let $\underline{X}$ be a del Pezzo surface of degree $5$ and $C$ is an irreducible curve of anticanonical degree $10$. Then $C^2 = 8, 10, 12, 14, 16, 18$ or $20$.
    When $C^2 = 8$, there exists a birational morphism $\pi: \underline{X} \to \underline{X}'$ to $\PP^1 \times \PP^1$ such that $C$ is a member of pullback of a member in $(4,1)$ or $(1,4)$-class in $\PP^1 \times \PP^1$. 
    When $C^2 = 10$, no such effective curve exists.
    When $C^2 = 12$, there exists a birational morphism $\pi: \underline{X} \to \underline{X}'$ to blow up of $\PP^2$ at one point such that $C$ is a member of $|4\pi^*H' - 2\pi^*E|$, where $H'$ is the pullback of a hyperplane class in $\PP^2$ via a blow up $\beta':\underline{X}' \to \PP^2$ at one point and $E$ is the unique exceptional divisor of $\beta'$.
    When $C^2 = 14$, there exists a birational morphism $\pi: \underline{X} \to \underline{X}'$ to blow up of $\PP^2$ at three points such that $C$ is a member of $|-\pi^*K_{\underline{X}'} + 2\pi^*F|$, for some smooth conic $F$.
    When $C^2 = 16$, then $C$ is either of the form $5H - 2E_i - 2E_j - E_k,$ $6H - 3E_i - 3E_j - E_k - E_l$ or $7H - 4E_i - 3E_j - 2E_k - 2E_l$, where $H$ is the pullback of a hyperplane class in $\PP^2$ by $\beta$, and $1 \leq i, j, k, l \leq 4$ are distinct indices.
    When $C^2 = 18$, there exists a birational morphism $\pi: \underline{X} \to \underline{X}'$ blow up of $\PP^2$ at one point such that $C$ is a member of $|-K_{\underline{X}} + 2\pi^*H' - \pi^*E|$, where $H'$ is the pullback of a hyperplane class in $\PP^2$ via a blow up $\beta':\underline{X}' \to \PP^2$ at one point and $E$ is the unique exceptional divisor of $\beta'$.
    When $C^2 = 20$, then $C \in |-2K_{\underline{X}}|$.
\end{proposition}

\begin{proof}
    \begin{enumerate}
        \item \textit{Case 1: $C^2 = 8$.} In this case $2C$ is big and nef. Hence using Kawamata-Viehweg vanishing and Riemann-Roch we can conclude that $2C + K_{\underline{X}}$ is effective. We also have $$C \cdot (2C + K_{\underline{X}}) = 6, \quad -K_{\underline{X}} \cdot (2C + K_{\underline{X}}) = 15, \quad (2C + K_{\underline{X}})^2 = -3.$$
        After writing the Zariski decomposition of the divisor $2C + K_{\underline{X}}$, by similar calculations as in proposition \ref{dp3d5} we can conclude that the number of $N_i$'s is $3$. Let $\pi: \underline{X} \to \underline{X}'$ be a birational contraction of $N_i$'s.
        This is only possible when $\underline{X}'$ is blow-up of $\PP^2$ at one point or is $\PP^1 \times \PP^1$. Let $H'$ be the pullback of a hyperplane class in $\PP^2$ via a blow up $\beta':\underline{X}' \to \PP^2$ at one point and $E$ be the unique exceptional divisor of $\beta'$, then $2C' \sim 9H' - 7E$, but the divisor $9H' - 7E$ is not divisible by $2$. Hence this case is not possible. Hence $\underline{X}'$ must be $\PP^1 \times \PP^1$ and $C'$ is either a member of the $(4,1)$ or $(1,4)$-class in $\PP^1 \times \PP^1$. 
        For different contraction maps, they are one of the following forms. Using \ref{prog} we get:
        \begin{enumerate}
            \item $5H - 4E_i - E_j$, where $1 \leq i, j \leq 4, i \neq j$.
            \item $6H - 5E_i - E_j - E_k - E_l$.
            \item $9H - 5E_i - 4E_j - 4E_k - 4E_l$, where $1 \leq i, j, k, l \leq 4$ are distinct indices.
        \end{enumerate}
        \item \label{dp5d5_10} \textit{Case 2: $C^2 = 10$.} Using \ref{prog} we get no solution for such curves which are effective. 
        \item \textit{Case 3: $C^2 = 12$.} In this case $C$ is big and nef. Using Kawamata-Viehweg vanishing and Riemann-Roch we can conclude that $|C + K_{\underline{X}}|$ defines a pencil. We also have $$C \cdot (C + K_{\underline{X}}) = 2, \quad -K_{\underline{X}} \cdot (C + K_{\underline{X}}) = 5, \quad (C + K_{\underline{X}})^2 = -3.$$
        After writing the Zariski decomposition of the divisor $C + K_{\underline{X}}$, by similar calculations as in proposition \ref{dp3d5} we can conclude that the number of $N_i$'s is $3$. Let $\pi: \underline{X} \to \underline{X}'$ be a birational contraction of $N_i$'s.
        This is only possible when $\underline{X}'$ is blow-up of $\PP^2$ at one point. Let $H'$ be the pullback of a hyperplane class in $\PP^2$ via a blow up $\beta':\underline{X}' \to \PP^2$ at one point and $E$ be the unique exceptional divisor of $\beta'$, then $C' \sim 4H' - 2E$ and $C \in |4\pi^*H' - 2\pi^*E|$.
        For different contraction maps, they are one of the following forms. Using \ref{prog} we get:
        \begin{enumerate}
            \item $4H - 2E_i$, where $1 \leq i \leq 4$.
            \item $5H - 3E_i - 2E_j$, where $1 \leq i, j \leq 4, i \neq j$.
            \item $6H - 4E_i - 2E_j - 2E_k$, where $1 \leq i, j, k \leq 4$ are distinct indices.
            \item $7H - 5E_i - 2E_j - 2E_k - 2E_l$.
            \item $8H - 4E_i - 4E_j - 4E_k - 2E_l;$ $8H - 5E_i - 3E_j - 3E_k - 3E_l$, where $1 \leq i, j, k, l \leq 4$ are distinct indices.
        \end{enumerate}
        \item \textit{Case 4: $C^2 = 14$.} In this case $C$ is big and nef. Using Kawamata-Viehweg vanishing and Riemann-Roch we can conclude that $C + K_{\underline{X}}$ is effective. We also have $$C \cdot (C + K_{\underline{X}}) = 4, \quad -K_{\underline{X}} \cdot (C + K_{\underline{X}}) = 5, \quad (C + K_{\underline{X}})^2 = -1.$$
        After writing the Zariski decomposition of the divisor $C + K_{\underline{X}}$, by similar calculations as in proposition \ref{dp3d5} we can conclude that the number of $N_i$'s is $1$. Let $\pi: \underline{X} \to \underline{X}'$ be a birational contraction of $N_i$ by $C$.
        This is only possible when $\underline{X}'$ is blow-up of $\PP^2$ at three points. Let $C'$ be the pushforward of $C$ in $\underline{X}'$. In $\underline{X}'$ we have $(C' + K_{\underline{X}'})^2 = 0$ and $-K_{\underline{X}'} \cdot (C' + K_{\underline{X}'}) = 4$. Using similar discussions as in proposition \ref{dp4d8}, case 4, we can conclude that $C'$ is linearly equivalent to $-K_{\underline{X}'} + 2F$ for some smooth conic $F$. Hence $C$ is a member of $|-\pi^*K_{\underline{X}'} + 2\pi^*F|$ for some smooth conic $F$.
        They are one of the following forms. Using \ref{prog} we get:
        \begin{enumerate}
            \item $4H - E_i - E_j$, where $1 \leq i, j \leq 4, i \neq j$.
            \item $5H - 3E_i - E_j - E_k$, where $1 \leq i, j, k \leq 4$ are distinct indices.
            \item $6H - 3E_i - 3E_j - 2E_k;$ $6H - 4E_i - 2E_j - E_k - E_l$.
            \item $7H - 4E_i - 3E_j - 3E_k - E_l$.
            \item $8H - 4E_i - 4E_j - 3E_k - 3E_l$, where $1 \leq i, j, k, l \leq 4$ are distinct indices.
        \end{enumerate}
        \item \label{dp5d5_16} \textit{Case 5: $C^2 = 16$.} Using \ref{prog} we get the following forms of curves of this type.
        \begin{enumerate}
            \item $5H - 2E_i - 2E_j - E_k$, where $1 \leq i, j, k \leq 4$ are distinct indices.
            \item $6H - 3E_i - 3E_j - E_k - E_l$.
            \item $7H - 4E_i - 3E_j - 2E_k - 2E_l$, where $1 \leq i, j, k, l \leq 4$ are distinct indices.
        \end{enumerate}
        \item \textit{Case 6: $C^2 = 18$.} In this case $(C+K_{\underline{X}})^2 = 3$ and $-K_{\underline{X}} \cdot (C + K_{\underline{X}}) = 5$. Then there exists a birational morphism $\pi: \underline{X} \to \underline{X}'$ to blow up of $\PP^2$ at one point such that $C$ is a member of $|-K_{\underline{X}} + 2\pi^*H' - \pi^*E|,$ where $H'$ is the pullback of a hyperplane class in $\PP^2$ via a blow up $\beta':\underline{X}' \to \PP^2$ at one point and $E$ is the exceptional divisor of $\beta'$. Then we have the following forms of $C$. Using \ref{prog} we get:
        \begin{enumerate}
            \item $5H - 2E_i - E_j - E_k - E_l$.
            \item $6H - 3E_i - 2E_j - 2E_k - E_l$.
            \item $7H - 3E_i - 3E_j - 3E_k - 2E_l$, where $1 \leq i, j, k, l \leq 4$ are distinct indices.
        \end{enumerate}
        \item \textit{Case 7: $C^2 = 20$.} In this case $C \in |-2K_{\underline{X}}|.$
    \end{enumerate}
\end{proof}

\subsection{Classification of irreducible boundaries.} Let $\underline{X}$ be a del Pezzo surface of degree $5$ and $\Delta$ be a $\Q$-divisor. We classify different positivity properties of $\Delta = \epsilon D$ when $D$ is an irreducible curve of anticanonical degree at most $10$ and $\epsilon = 1 - 1/m.$

Using Lemma \ref{thelemma} for a del Pezzo Surface of degree $5$, we can conclude that a divisor $L$ is ample (nef) if and only if $L \cdot E > 0$ (resp. $L \cdot E \geq 0$) for any $(-1)$-curve $E$. 
We use \ref{prog} to calculate the intersection numbers.

\begin{itemize}
    \item \textit{Case of $-K_{\underline{X}} \cdot D = 1$.} We assume that $D$ is a $(-1)$-curve. In this case $D^2 = -1$. Since $D \cdot E = 0$ or $1$ for any $(-1)$-curve $E$ distinct from $D$, we conclude that $-(K_{\underline{X}} + \epsilon D)$ is ample for any $\epsilon > 0.$ In this case $-(K_{\underline{X}} + D)$ is big and nef but not ample.
    \item \textit{Case of $-K_{\underline{X}} \cdot D = 2$.} We assume that $D$ is a smooth fiber of a conic fibration. In this case, $D \cdot E = 0$ or $1$ for any $(-1)$-curve $E$. Therefore we conclude that $-(K_{\underline{X}} + \epsilon D)$ is ample for any $\epsilon > 0.$
    \item \textit{Case of $-K_{\underline{X}} \cdot D = 3$.} We assume that $D$ is pullback of a line in $\PP^2$. In this case $D \cdot E = 0$ or $1$ for any $(-1)$-curve $E$. Therefore we conclude that $-(K_{\underline{X}} + \epsilon D)$ is ample for any $\epsilon > 0.$
    \item \textit{Case of $-K_{\underline{X}} \cdot D = 4$.} We assume that $D$ is pullback of the $(1,1)$-class in $\PP^1 \times \PP^1$. Then $D \cdot E = 0, 1$ or $2$ for any $(-1)$-curve $E$. Therefore we conclude that $-(K_{\underline{X}} + \epsilon D)$ is big and nef only for $\epsilon = \frac{1}{2}.$
    \item \textit{Case of $-K_{\underline{X}} \cdot D = 5$.} First we assume that $D^2 = 3$. Then $D$ is pullback of a divisor of the form $2H' - E$ as described above. In this case $D \cdot F = 0, 1$ or $2$ for any $(-1)$-curve $F$. Hence we conclude, $-(K_{\underline{X}} + \epsilon D)$ is big and nef only for $\epsilon = \frac{1}{2}.$

    Next we assume that $D$ is a member of the anticanonical class. Since $D \cdot E = 1$ for any $(-1)$-curve $E$, we conclude that $-(K_{\underline{X}} + \epsilon D)$ is ample for all $\epsilon > 0.$
    \item \textit{Case of $-K_{\underline{X}} \cdot D = 6.$} We first assume that $D^2 = 4$. In this case $D \cdot E = 0$ or $2$ for some curves and so $-(K_{\underline{X}} + \epsilon D)$ is nef for $\epsilon = \frac{1}{2}$, but it is not big, and $D \cdot F = 0, 1, 2$ or $3$ for the other type of curves, for all $(-1)$-curves $F$ and so $-(K_{\underline{X}} + \epsilon D)$ is not nef for any $\epsilon > 0$.

    Next we assume that $D$ is a member of pullback of the anticanonical class from a smooth sextic surface. In this case $D \cdot E = 0, 1$ or $2$ for any $(-1)$-curve $E$. Hence $-(K_{\underline{X}} + \epsilon D)$ is big and nef only for $\epsilon = \frac{1}{2}.$
    \item \textit{Case of $-K_{\underline{X}} \cdot D = 7$.} We first assume that $D^2 = 5.$ Then $D$ is pullback of a divisor of the form $3H' - 2E$ as described above. In this case $D \cdot F = 0, 1, 2$ or $3$ for any $(-1)$-curve $F$. Hence $-(K_{\underline{X}} + \epsilon D)$ is not nef for any $\epsilon > 0.$

    Next we assume that $D^2 = 7$. In this case $D$ is pullback of a member in the anticanonical class of a del Pezzo surface of degree $7$. Then $D \cdot E = 0, 1 ,2$ or $3$ for any $(-1)$-curve $E$. So $-(K_{\underline{X}} + \epsilon D)$ is not nef for any $\epsilon > 0.$

    Next we assume that $D^2 = 9$. In this case $D + K_{\underline{X}}$ is linearly equivalent to a smooth conic. Then $D \cdot E = 1$ or $2$. Hence $-(K_{\underline{X}} + \epsilon D)$ is big and nef only for $\epsilon = \frac{1}{2}.$
    \item \textit{Case of $-K_{\underline{X}} \cdot D = 8$.} We first assume that $D^2 = 6$. Then $D$ is pullback of a divisor in the $(3,1)$ or $(1,3)$-class in $\PP^1 \times \PP^1$ as described above. Then $D \cdot F = 0, 1, 3$ or $4$ for any $(-1)$-curve $F$. Hence we conclude that $-(K_{\underline{X}} + \epsilon D)$ is not nef for any $\epsilon > 0.$

    Next we assume that $D^2 = 8$. Then $D$ is pullback of a member in the anticanonical class of blow up of $\PP^2$ at one point. In this case, $D \cdot E = 0, 2$ or $4$ for any $(-1)$-curve $E$, for two type of curves, and $D \cdot E = 0, 1, 2$ or $3$ for other types, for any $(-1)$-curve $E$. In any case, we conclude that $-(K_{\underline{X}} + \epsilon D)$ is not nef for any $\epsilon > 0.$

    Next we assume that $D^2 = 10$. In this case $D \cdot F = 0, 1, 2$ or $3$ for any $(-1)$-curve $F$. Hence we conclude that $-(K_{\underline{X}} + \epsilon D)$ is not nef for any $\epsilon > 0.$

    Next we assume that $D^2 = 12$. Then $D + K_{\underline{X}}$ is pullback of a line in $\PP^2$ as described above. In this case $D \cdot E = 1$ or $2$ for any $(-1)$-curve $E$. Hence we conclude that $-(K_{\underline{X}} + \epsilon D)$ is nef only for $\epsilon = \frac{1}{2},$ but it is not big.
    \item \textit{Case of $-K_{\underline{X}} \cdot D = 9$.} We first assume that $D^2 = 7$. Then $D$ is pullback of a divisor of the form $4H - 3E$ as described above, where $H$ is the pullback of a hyperplane class in $\PP^2$  and $E$ is the exceptional divisor. In this case $D \cdot E = 0, 1, 3$ or $4$. Hence we conclude that $-(K_{\underline{X}} + \epsilon D)$ is not nef for any $\epsilon > 0.$

    Next we assume that $D^2 = 9.$ Then $D$ is pullback of a member in the anticanonical class of $\PP^2$. In this case $D \cdot E = 0$ or $3$ for any $(-1)$-curve $E$. Hence we conclude that $-(K_{\underline{X}} + \epsilon D)$ is not nef for any $\epsilon > 0.$

    Next we assume that $D^2 = 11$. In this case $D \cdot F = 0, 1, 2, 3$ or $4$ for any $(-1)$-curve $E$. So we conclude that $-(K_{\underline{X}} + \epsilon D)$ is not nef for any $\epsilon > 0.$

    Next we assume that $D^2 = 13$. Then $D + K_{\underline{X}}$ is linearly equivalent to $2F$ where $F$ is a smooth conic. In this case $D \cdot E = 1$ or $3$ for all $(-1)$-curve $E$, for two type of conics and $D \cdot E = 0, 1, 2$ or $3$ for all $(-1)$-curve $E$ for other type of conics. In any case, we conclude that $-(K_{\underline{X}} + \epsilon D)$ is not nef for any $\epsilon > 0.$

    Next we assume that $D^2 = 15.$ Then $D + K_{\underline{X}}$ is pullback of the $(1,1)$-class in $\PP^1 \times \PP^1$. In this case $D \cdot E = 1, 2$ or $3$ for any $(-1)$-curve $E$. Hence $-(K_{\underline{X}}+\epsilon D)$ is not nef for any $\epsilon > 0.$
    \item \textit{Case of $-K_{\underline{X}} \cdot D = 10$.} We first assume that $D^2 = 8$. Then $D$ is pullback of a divisor of the form $(4,1)$ or $(1,4)$-class in $\PP^1 \times \PP^1$. In this case $D \cdot F$ can be as high as $5$, hence $-(K_{\underline{X}} + \epsilon D)$ can not be nef for any $\epsilon > 0.$

    Next we assume that $D^2 = 12.$ Then $D$ is pullback of a divisor of the form $4H' - 2E$ as explained above, where $H'$ is the pullback of a hyperplane class in $\PP^2$ via a blow up $\beta':\underline{X}' \to \PP^2$ at one point and $E$ is the exceptional curve of $\beta'$. In this case $D \cdot F$ can be as high as $5$, hence $-(K_{\underline{X}} + \epsilon D)$ can not be nef for any $\epsilon > 0.$

    Next we assume that $D^2 = 14.$ In this case $D \cdot E = 0, 1, 2, 3$ or $4$ for any $(-1)$-curve $E$. Hence $-(K_{\underline{X}} + \epsilon D)$ is not nef for any $\epsilon > 0.$

    Next we assume that $D^2 = 16$. Then $D$ is one of the curves of the three types mentioned above. In this case $D \cdot E = 0, 1, 2, 3$ or $4$ for any $(-1)$-curve $E$. Hence $-(K_{\underline{X}} + \epsilon D)$ is not nef for any $\epsilon > 0.$

    Next we assume that $D^2 = 18.$ In this case $D + K_{\underline{X}}$ is pullback of a divisor of the form $2H' - E$ as explained above, where $H'$ is the pullback of a hyperplane class in $\PP^2$ via a blow up $\beta':\underline{X}' \to \PP^2$ at one point and $E$ is exceptional curve of $\beta'$. Then $D \cdot F = 1, 2$ or $3$ for any $(-1)$-curve $F$. Hence $-(K_{\underline{X}} + \epsilon D)$ is not nef for any $\epsilon > 0.$

    Next we assume that $D \in |-2K_{\underline{X}}|$. Since $D \cdot E = 2$ for all $(-1)$-curves $E$, we conclude that $-(K_{\underline{X}} + \epsilon D)$ is nef only for $\epsilon = \frac{1}{2}$, but it is not big.
\end{itemize}

In summary, we conclude:

\begin{proposition}
    Let $(X, \epsilon D)$ be a del Pezzo orbifold whose underlying surface $\underline{X}$ is a del Pezzo surface of degree $5$ and $D$ is an irreducible curve of degree less than or equal to $10$, and $\epsilon = 1 - 1/m, m \in \Z_{\geq 2}$. Then $-(K_{\underline{X}} + \epsilon D)$ is ample if and only if $D$ is one of the following.
    \begin{itemize}
        \item a $(-1)$-curve and $\epsilon > 0.$
        \item a smooth fiber of a conic fibration and $\epsilon > 0.$
        \item pullback of a line in $\PP^2$ via a birational morphism and $\epsilon > 0.$
        \item member of the anticanonical class and $\epsilon > 0.$
    \end{itemize}
    Also, $-(K_{\underline{X}} + \epsilon D)$ is nef when $\epsilon = \frac{1}{2}$ and $D$ is one of the following.
    \begin{itemize}
        \item pullback of the $(1,1)$-class in $\PP^1 \times \PP^1$. In this case $-(K_{\underline{X}} + \epsilon D)$ is big.
        \item pullback of a divisor of the form $2H' - E$ via a birational morphism to blow up of $\PP^2$ at one point, where $H'$ is the pullback of a hyperplane class in $\PP^2$ via a blow up $\beta'$ as mentioned above and $E$ is the exceptional divisor of $\beta'$. In this case $-(K_{\underline{X}} + \epsilon D)$ is big.
        \item pullback of a divisor in the $(1,2)$ or $(2,1)$-class in $\PP^1 \times \PP^1$ via a birational morphism to $\PP^1 \times \PP^1$. In particular, they are of the form:
        $2H$ or $4H - 2E_i - 2E_j - 2E_k$ where $H$ is the pullback of a hyperplane class in $\PP^2$ via $\beta$. In this case $-(K_{\underline{X}} + \epsilon D)$ is not big.
        \item member of pullback of the anticanonical class of a smooth sextic del Pezzo surface. In this case $-(K_{\underline{X}} + \epsilon D)$ is big.
        \item $D$ is linearly equivalent to $-K_{\underline{X}} + F$, for a smooth conic $F$. In this case $-(K_{\underline{X}} + \epsilon D)$ is big.
        \item $D + K_{\underline{X}}$ is pullback of a line in $\PP^2$. In this case $-(K_{\underline{X}} + \epsilon D)$ is not big.
        \item member of $|-2K_{\underline{X}}|$. In this case $-(K_{\underline{X}} + \epsilon D)$ is not big.
    \end{itemize}
\end{proposition}

\begin{landscape}
    {\scriptsize
        \renewcommand{\arraystretch}{2}
        \begin{longtable}{|c|c|c|c|c|c|c|}
        \caption{Low degree irreducible curves on del Pezzo surfaces}\\
            \hline
            \multirow{1}{*}{\makecell{\textbf{anticanonical} \\ \textbf{degree}}}
            & \makecell{\textbf{self intersection} \\ \textbf{number}} & \hyperref[dp1]{\makecell{\textbf{del Pezzo surface} \\ \textbf{of degree} $1$}} & \hyperref[dp2]{\makecell{\textbf{del Pezzo surface} \\ \textbf{of degree} $2$}} & \hyperref[dp3]{\makecell{\textbf{del Pezzo surface} \\ \textbf{of degree} $3$}} & \hyperref[dp4]{\makecell{\textbf{del Pezzo surface} \\ \textbf{of degree} $4$}} & \hyperref[dp5]{\makecell{\textbf{del Pezzo surface} \\ \textbf{of degree} $5$}} \\ \hline
            
            \multirow{2}{*}{$-K_{\underline{X}} \cdot C = 1$}
            & $C^2 = -1$ & $(-1)$-curve & $(-1)$-curve & $\textcolor{blue}{(-1)}$\textcolor{blue}{-curve} & $\textcolor{blue}{(-1)}$\textcolor{blue}{-curve} & $\textcolor{blue}{(-1)}$\textcolor{blue}{-curve} \\ \cline{2-7}
            & $C^2 = 1$ & $\textcolor{blue}{C \in \lvert-K_{\underline{X}}\rvert}$ &  &  &  &  \\ \hline
            
            \multirow{3}{*}{$-K_{\underline{X}} \cdot C = 2$}
            & $C^2 = 0$ & $F$ (conic) & $F$ (conic) & $F$ (conic) & $\textcolor{blue}{F \text{ (conic)}}$ & $\textcolor{blue}{F \text{ (conic)}}$ \\ \cline{2-7}
            & $C^2 = 2$ & \makecell{\vspace{0.4mm} \\ $C \in \lvert -\pi^* K_{\underline{X}'} \rvert$ \\ $\underline{X}' = $ del Pezzo surface \\ of degree $2$ \\ \vspace{0.4mm}} & $\textcolor{blue}{C \in \lvert -K_{\underline{X}} \rvert}$ &  &  &  \\ \cline{2-7} 
            & $C^2 = 4$ & $C \in \lvert -2K_{\underline{X}} \rvert$ &  &  &  &  \\ \hline
            
            \multirow{2}{*}{$-K_{\underline{X}} \cdot C = 3$}
            & $C^2 = 1$ & & \makecell{\vspace{0.4mm} \\ $C \in \lvert \pi^*H \rvert$ \\ $\pi: \underline{X} \to \mathbb{P}^2$ \\ \vspace{0.4mm}} & \makecell{\vspace{0.4mm} \\ $C \in \lvert \pi^*H \rvert$ \\ $\pi: \underline{X} \to \mathbb{P}^2$ \\ \vspace{0.4mm}} & \makecell{\vspace{0.4mm} \\ $C \in \lvert \pi^*H \rvert$ \\ $\pi: \underline{X} \to \mathbb{P}^2$ \\ \vspace{0.4mm}} & \makecell{\vspace{0.4mm} \\ $\textcolor{blue}{C \in \lvert \pi^*H \rvert}$ \\ $\textcolor{blue}{\pi: \underline{X} \to \mathbb{P}^2}$ \\ \vspace{0.4mm}}  \\ \cline{2-2}\cline{4-7}
            & $C^2 = 3$ & & \makecell{\vspace{0.4mm} \\ $C \in \lvert -\pi^* K_{\underline{X}'} \rvert$ \\ $\underline{X}' = $ del Pezzo surface \\ of degree $3$ \\ \vspace{0.4mm}} & $\textcolor{blue}{C \in \lvert -K_{\underline{X}} \rvert}$ &  &  \\ \cline{2-7} \hline
            
            \multirow{4}{*}{$-K_{\underline{X}} \cdot C = 4$}
            & $C^2 = 2$ &  & \makecell{\vspace{0.4mm} \\ $C \in \lvert \pi^*(1,1) \rvert$ \\ $\pi: \underline{X} \to \mathbb{P}^1 \times \mathbb{P}^1$ \\ \vspace{0.4mm}} & \makecell{\vspace{0.4mm} \\ $C \in \lvert \pi^*(1,1) \rvert$ \\ $\pi: \underline{X} \to \mathbb{P}^1 \times \mathbb{P}^1$ \\ \vspace{0.4mm}} & \makecell{\vspace{0.4mm} \\ $C \in \lvert \pi^*(1,1) \rvert$ \\ $\pi: \underline{X} \to \mathbb{P}^1 \times \mathbb{P}^1$ \\ \vspace{0.4mm}} & \makecell{\vspace{0.4mm} \\ $C \in \lvert \pi^*(1,1) \rvert$ \\ $\pi: \underline{X} \to \mathbb{P}^1 \times \mathbb{P}^1$ \\ \vspace{0.4mm}}  \\ \cline{2-2} \cline{4-7}
            & $C^2 = 4$ &  & \makecell{\vspace{0.4mm} \\ $C \in \lvert -\pi^*K_{\underline{X}'} \rvert$ \\ $\underline{X}' = $ del Pezzo surface \\ of degree $4$ \\ \vspace{0.4mm}} & \makecell{\vspace{0.4mm} \\ $C \in \lvert -\pi^*K_{\underline{X}'} \rvert$ \\ $\underline{X}' = $ del Pezzo surface \\ of degree $4$ \\ \vspace{0.4mm}} & $\textcolor{blue}{C \in \lvert -K_{\underline{X}} \rvert}$ &   \\ \cline{2-2} \cline{4-7}
            & $C^2 = 6$ &  & $C \sim -K_{\underline{X}} + F$ &  &  &  \\ \cline{2-2} \cline{4-4}
            & $C^2 = 8$ &  & $C \in \lvert -2K_{\underline{X}} \rvert$ &  &  &  \\ \hline
        \end{longtable}
        \newpage
        \begin{longtable}{|c|c|c|c|c|c|c|}
            \hline
            \multirow{1}{*}{\makecell{\textbf{anticanonical} \\ \textbf{degree}}}
            & \makecell{\textbf{self intersection} \\ \textbf{number}} & \hyperref[dp1]{\makecell{\textbf{del Pezzo surface} \\ \textbf{of degree} $1$}} & \hyperref[dp2]{\makecell{\textbf{del Pezzo surface} \\ \textbf{of degree} $2$}} & \hyperref[dp3]{\makecell{\textbf{del Pezzo surface} \\ \textbf{of degree} $3$}} & \hyperref[dp4]{\makecell{\textbf{del Pezzo surface} \\ \textbf{of degree} $4$}} & \hyperref[dp5]{\makecell{\textbf{del Pezzo surface} \\ \textbf{of degree} $5$}} \\ \hline
            
            \multirow{3}{*}{$-K_{\underline{X}} \cdot C = 5$}
            & $C^2 = 3$ &  &  & \makecell{\vspace{0.4mm} \\ $C \in \lvert 2\pi^*H' - \pi^*E \rvert$ \\ $\pi: \underline{X} \to \underline{X}'$ \\ $\underline{X}' = \mathbb{P}^2$ blown up \\ at one point \\ \vspace{0.4mm}} & \makecell{\vspace{0.4mm} \\ $C \in \lvert 2\pi^*H' - \pi^*E \rvert$ \\ $\pi: \underline{X} \to \underline{X}'$ \\ $\underline{X}' = \mathbb{P}^2$ blown up \\ at one point  \\ \vspace{0.4mm}} & \makecell{\vspace{0.4mm} \\ $C \in \lvert 2\pi^*H' - \pi^*E \rvert$ \\ $\pi: \underline{X} \to \underline{X}'$ \\ $\underline{X}' = \mathbb{P}^2$ blown up \\ at one point \\ \vspace{0.4mm}} \\ \cline{2-2} \cline{5-7}
            & $C^2 = 5$ &  &  & \makecell{\vspace{0.4mm} \\ $C \in \lvert -\pi^*K_{\underline{X}'} \rvert$ \\ $\underline{X}' = $ del Pezzo surface \\ of degree $5$ \\ \vspace{0.4mm}} & \makecell{\vspace{0.4mm} \\ $C \in \lvert -\pi^*K_{\underline{X}'} \rvert$ \\ $\underline{X}' = $ del Pezzo surface \\ of degree $5$ \\ \vspace{0.4mm}} & $\textcolor{blue}{C \in \lvert -K_{\underline{X}} \rvert}$ \\ \cline{2-2} \cline{5-7}
            & $C^2 = 7$ &  &  & $C \sim -K_{\underline{X}} + F$ &  &  \\ \hline
            
            \multirow{5}{*}{$-K_{\underline{X}} \cdot C = 6$}
            & $C^2 = 4$ &  &  & \makecell{\vspace{0.4mm} \\ $C \in \lvert \pi^*(1,2) \rvert$  \\ or $\lvert \pi^*(2,1) \rvert$ \\ $\pi: \underline{X} \to \mathbb{P}^1 \times \mathbb{P}^1$ \\ \vspace{0.4mm}} & \makecell{\vspace{0.4mm} \\ $C \in \lvert \pi^*(1,2) \rvert$  \\ or $\lvert \pi^*(2,1) \rvert$ \\ $\pi: \underline{X} \to \mathbb{P}^1 \times \mathbb{P}^1$ \\ \vspace{0.4mm}} & \makecell{\vspace{0.4mm} \\ $C \in \lvert \pi^*(1,2) \rvert$ \\ or $\lvert \pi^*(2,1) \rvert$ \\ $\pi: \underline{X} \to \mathbb{P}^1 \times \mathbb{P}^1$ \\ \vspace{0.4mm}} \\ \cline{2-2} \cline{5-7}
            & $C^2 = 6$ &  &  & \makecell{\vspace{0.4mm} \\ $C \in \lvert -\pi^*K_{\underline{X}'} \rvert$ \\ $\underline{X}' = $ del Pezzo surface  \\ of degree $6$ \\ \vspace{0.4mm}} & \makecell{\vspace{0.4mm} \\ $C \in \lvert -\pi^*K_{\underline{X}'} \rvert$ \\ $\underline{X}' = $ del Pezzo surface \\ of degree $6$ \\ \vspace{0.4mm}} &  \makecell{\vspace{0.4mm} \\ $C \in \lvert -\pi^*K_{\underline{X}'} \rvert$ \\ $\underline{X}' = $ del Pezzo surface \\ of degree $6$ \\ \vspace{0.4mm}} \\ \cline{2-2} \cline{5-7}
            & $C^2 = 8$ &  &  & \makecell{\vspace{0.4mm} \\ $C \in \lvert -\pi^*K_{\underline{X}'} + \pi^*F \rvert$ \\ $\underline{X}' = $ del Pezzo surface \\ of degree $4$ \\ \vspace{0.4mm}} & $C \sim -K_{\underline{X}} + F$ &  \\ \cline{2-2} \cline{5-6}
            & $C^2 = 10$ &  &  & \makecell{\vspace{0.4mm} \\ $C \in \lvert -K_{\underline{X}} + \pi^*H' \rvert$ \\ $\pi: \underline{X} \to \mathbb{P}^2$ \\ \vspace{0.4mm}} &  &  \\ \cline{2-2} \cline{5-5}
            & $C^2 = 12$ &  &  & $C \in \lvert -2K_{\underline{X}} \rvert$ &  &  \\ \hline
        \end{longtable}
        \newpage
        \begin{longtable}{|c|c|c|c|c|c|c|}
            \hline
            \multirow{1}{*}{\makecell{\textbf{anticanonical} \\ \textbf{degree}}}
            & \makecell{\textbf{self intersection} \\ \textbf{number}} & \hyperref[dp1]{\makecell{\textbf{del Pezzo surface} \\ \textbf{of degree} $1$}} & \hyperref[dp2]{\makecell{\textbf{del Pezzo surface} \\ \textbf{of degree} $2$}} & \hyperref[dp3]{\makecell{\textbf{del Pezzo surface} \\ \textbf{of degree} $3$}} & \hyperref[dp4]{\makecell{\textbf{del Pezzo surface} \\ \textbf{of degree} $4$}} & \hyperref[dp5]{\makecell{\textbf{del Pezzo surface} \\ \textbf{of degree} $5$}} \\ \hline
            
            \multirow{4}{*}{$-K_{\underline{X}} \cdot C = 7$}
            & $C^2 = 5$ &  &  &  & \makecell{\vspace{0.4mm} \\ $C \in \lvert 3\pi^*H' - 2\pi^*E \rvert$ \\ $\pi: \underline{X} \to \underline{X}'$ \\ $\underline{X}' = \mathbb{P}^2$ blown up \\ at one point \\ \vspace{0.4mm}} & \makecell{\vspace{0.4mm} \\ $C \in \lvert 3\pi^*H' - 2\pi^*E \rvert$ \\ $\pi: \underline{X} \to \underline{X}'$ \\ $\underline{X}' = \mathbb{P}^2$ blown up \\ at one point \\ \vspace{0.4mm}} \\ \cline{2-2} \cline{6-7}
            & $C^2 = 7$ &  &  &  & \makecell{\vspace{0.4mm} \\ $C \in \lvert -\pi^*K_{\underline{X}'} \rvert$ \\ $\underline{X}' = $ del Pezzo surface \\ of degree $7$ \\ \vspace{0.4mm}} & \makecell{\vspace{0.4mm} \\ $C \in \lvert -\pi^*K_{\underline{X}'} \rvert$ \\ $\underline{X}' = $ del Pezzo surface \\ of degree $7$ \\ \vspace{0.4mm}} \\ \cline{2-2} \cline{6-7}
            & $C^2 = 9$ &  &  &  & \makecell{\vspace{0.4mm} \\ $C \in \lvert -\pi^*K_{\underline{X}'} + \pi^*F \rvert$ \\ $\underline{X}' = $ del Pezzo surface \\ of degree $5$ \\ \vspace{0.4mm}} & $C \sim -K_{\underline{X}} + F$ \\ \cline{2-2} \cline{6-7}
            & $C^2 = 11$ &  &  &  & \makecell{\vspace{0.4mm} \\ $C \in \lvert -K_{\underline{X}} + \pi^*H' \rvert$ \\ $\pi: \underline{X} \to \mathbb{P}^2$ \\ \vspace{0.4mm}} &  \\ \hline
            
            \multirow{6}{*}{$-K_{\underline{X}} \cdot C = 8$}
            & $C^2 = 6$ &  &  &  & \makecell{\vspace{0.4mm} \\ $C \in \lvert \pi^*(1,3) \rvert$  \\ or $\lvert \pi^*(3,1) \rvert$ \\ $\pi: \underline{X} \to \mathbb{P}^1 \times \mathbb{P}^1$ \\ \vspace{0.4mm}} & \makecell{\vspace{0.4mm} \\ $C \in \lvert \pi^*(1,3) \rvert$  \\ or $\lvert \pi^*(3,1) \rvert$ \\ $\pi: \underline{X} \to \mathbb{P}^1 \times \mathbb{P}^1$ \\ \vspace{0.4mm}}  \\ \cline{2-2} \cline{6-7}
            & $C^2 = 8$ &  &  &  & \makecell{\vspace{0.4mm} \\ $C \in \lvert -\pi^*K_{\underline{X}'} \rvert$ \\ $\underline{X}' = $ del Pezzo surface \\ of degree $8$ \\ \vspace{0.4mm}} &  \makecell{\vspace{0.4mm} \\ $C \in \lvert -\pi^*K_{\underline{X}'} \rvert$ \\ $\underline{X}' = $ del Pezzo surface \\ of degree $8$ \\ \vspace{0.4mm}} \\ \cline{2-2} \cline{6-7}
            & $C^2 = 10$ &  &  &  & \makecell{\vspace{0.4mm} \\ $C \in \lvert -\pi^*K_{\underline{X}'} + \pi^*F \rvert$ \\ $\underline{X}' = $ del Pezzo surface \\ of degree $6$ \\ \vspace{0.4mm}} & \makecell{\vspace{0.4mm} \\ $C \in \lvert -\pi^*K_{\underline{X}'} + \pi^*F \rvert$ \\ $\underline{X}' = $ del Pezzo surface \\ of degree $6$ \\ \vspace{0.4mm}}  \\ \cline{2-2} \cline{6-7} \hline
        \end{longtable}
        \newpage
        \begin{longtable}{|c|c|c|c|c|c|c|}
            \hline
            \multirow{1}{*}{\makecell{\textbf{anticanonical} \\ \textbf{degree}}}
            & \makecell{\textbf{self intersection} \\ \textbf{number}} & \hyperref[dp1]{\makecell{\textbf{del Pezzo surface} \\ \textbf{of degree} $1$}} & \hyperref[dp2]{\makecell{\textbf{del Pezzo surface} \\ \textbf{of degree} $2$}} & \hyperref[dp3]{\makecell{\textbf{del Pezzo surface} \\ \textbf{of degree} $3$}} & \hyperref[dp4]{\makecell{\textbf{del Pezzo surface} \\ \textbf{of degree} $4$}} & \hyperref[dp5]{\makecell{\textbf{del Pezzo surface} \\ \textbf{of degree} $5$}} \\ \hline
            
            \multirow{2}{*}{$-K_{\underline{X}} \cdot C = 8$}
            & $C^2 = 12$ &  &  &  & $C \sim -K_{\underline{X}} + 2F$ & \makecell{\vspace{0.4mm} \\ $C \in \lvert -K_{\underline{X}} + \pi^*H' \rvert$ \\ $\pi: \underline{X} \to \mathbb{P}^2$ \\ \vspace{0.4mm}}  \\ \cline{2-2} \cline{6-7}
            & $C^2 = 14$ &  &  &  & \makecell{\vspace{0.4mm} \\ $C \in \lvert -K_{\underline{X}} + \pi^*(1,1) \rvert$ \\ $\pi: \underline{X} \to \mathbb{P}^1 \times \mathbb{P}^1$ \\ \vspace{0.4mm}} &  \\ \cline{2-2} \cline{6-7}
            & $C^2 = 16$ &  &  &  & $C \in \lvert-2K_{\underline{X}}\rvert$ & \\ \hline
        
            \multirow{5}{*}{$-K_{\underline{X}} \cdot C = 9$}
            & $C^2 = 7$ &  &  &  &  & \makecell{\vspace{0.4mm} \\ $C \in \lvert 4\pi^*H' - 3\pi^*E \rvert$ \\ $\pi: \underline{X} \to \underline{X}'$ \\ $\underline{X}' = \mathbb{P}^2$ blown up \\ at one point \\ \vspace{0.4mm}} \\ \cline{2-2} \cline{7-7}
            & $C^2 = 9$ &  &  &  &  & \makecell{\vspace{0.4mm} \\ $C \in \lvert -\pi^*K_{\underline{X}'} \rvert$ \\ $\underline{X}' = \mathbb{P}^2$ \\ \vspace{0.4mm}} \\ \cline{2-2} \cline{7-7}
            & $C^2 = 11$ &  &  &  &  & \makecell{\vspace{0.4mm} \\ $C \in \lvert -\pi^*K_{\underline{X}'} + \pi^*F \rvert$ \\ $\underline{X}' = $ del Pezzo surface \\ of degree $7$ \\ \vspace{0.4mm}} \\ \cline{2-2} \cline{7-7}
            & $C^2 = 13$ &  &  &  &  & $C \sim -K_{\underline{X}} + 2F$ \\ \cline{2-2} \cline{7-7}
            & $C^2 = 15$ &  &  &  &  & \makecell{\vspace{0.4mm} \\ $C \in \lvert -K_{\underline{X}} + \pi^*(1,1) \rvert$ \\ $\pi: \underline{X} \to \mathbb{P}^1 \times \mathbb{P}^1$ \\ \vspace{0.4mm}} \\ \hline
            
            \multirow{2}{*}{$-K_{\underline{X}} \cdot C = 10$}
            & $C^2 = 8$ &  &  &  &  & \makecell{\vspace{0.4mm} \\ $C \in \lvert \pi^*(1,4) \rvert$ \\ or $\lvert \pi^*(4,1) \rvert$ \\ $\pi: \underline{X} \to \mathbb{P}^1 \times \mathbb{P}^1$ \\ \vspace{0.4mm}}  \\ \cline{2-2} \cline{7-7}
            & $C^2 = 10$ &  &  &  &  & \hyperref[dp5d5_10]{No such curve exists} \\ \cline{2-2} \cline{7-7} \hline
        \end{longtable}
        \newpage
        \begin{longtable}{|c|c|c|c|c|c|c|}
            \hline
            \multirow{1}{*}{\makecell{\textbf{anticanonical} \\ \textbf{degree}}}
            & \makecell{\textbf{self intersection} \\ \textbf{number}} & \hyperref[dp1]{\makecell{\textbf{del Pezzo surface} \\ \textbf{of degree} $1$}} & \hyperref[dp2]{\makecell{\textbf{del Pezzo surface} \\ \textbf{of degree} $2$}} & \hyperref[dp3]{\makecell{\textbf{del Pezzo surface} \\ \textbf{of degree} $3$}} & \hyperref[dp4]{\makecell{\textbf{del Pezzo surface} \\ \textbf{of degree} $4$}} & \hyperref[dp5]{\makecell{\textbf{del Pezzo surface} \\ \textbf{of degree} $5$}} \\ \hline
            
            \multirow{5}{*}{$-K_{\underline{X}} \cdot C = 10$}
            & $C^2 = 12$ &  &  &  &  & \makecell{\vspace{0.4mm} \\ $C \in \lvert 4\pi^*H' - 2\pi^*E \rvert$ \\ $\pi: \underline{X} \to \underline{X}'$ \\ $\underline{X}' = \mathbb{P}^2$ blown up \\ at one point \\ \vspace{0.4mm}}  \\ \cline{2-2} \cline{7-7}
            & $C^2 = 14$ &  &  &  &  & \makecell{\vspace{0.4mm} \\ $C \in \lvert -\pi^*K_{\underline{X}'} + 2\pi^*F \rvert$ \\ $\underline{X}' = $ del Pezzo surface \\ of degree $6$ \\ \vspace{0.4mm}} \\ \cline{2-2} \cline{7-7}
            & $C^2 = 16$ &  &  &  &  & \hyperref[dp5d5_16]{explicit forms} \\ \cline{2-2} \cline{7-7}
            & $C^2 = 18$ &  &  &  &  & \makecell{\vspace{0.4mm} \\ $C \in \lvert -K_{\underline{X}} + \pi^*H' - \pi^*E \rvert$ \\ $\pi: \underline{X} \to \underline{X}'$ \\ $\underline{X}' = \mathbb{P}^2$ blown up  \\ at one point \\ \vspace{0.4mm}} \\ \cline{2-2} \cline{7-7}
            & $C^2 = 20$ &  &  &  &  & $C \in \lvert -2K_{\underline{X}} \rvert$ \\ \hline
            
        \end{longtable}
    }
\end{landscape}

\bibliographystyle{amsalpha}
\bibliography{ref}

\end{document}